\documentstyle{article}

 \font\tenmsb=msbm10 
 \font\sevenmsb=msbm7 
 \font\fivemsb=msbm5
 \newfam\msbfam
 \textfont\msbfam=\tenmsb
 \scriptfont\msbfam=\sevenmsb 
 \scriptscriptfont\msbfam=\fivemsb
 \def\Bbb{\fam\msbfam}

\font\teneuf=eufm10
\font\seveneuf=eufm7
\font\fiveeuf=eufm5
\newfam\euffam
\textfont\euffam=\teneuf
\scriptfont\euffam=\seveneuf
\scriptscriptfont\euffam=\fiveeuf

\title{On Distinguishing Quotients of Symmetric Groups}
\author{S.Shelah$^1$: Hebrew University of Jerusalem}
\date{and J. K. Truss: University of Leeds}  
\begin{document}
\maketitle  
\newtheorem{lemma}{Lemma}[section]
\newtheorem{theorem}[lemma]{Theorem}
\newtheorem{definition}[lemma]{Definition}
\newtheorem{remark}[lemma]{Remark}
\newtheorem{corollary}[lemma]{Corollary}
\setcounter{footnote}{1}\footnotetext{Research supported by the Israel Science Foundation,
administered by the Israel Academy of Sciences and Humanities, Publ. No. 605.}
\newcounter{number}

\begin{abstract}
A study is carried out of the elementary theory of quotients of symmetric groups in a similar spirit
to \cite{Shelah1}. Apart from the trivial and alternating subgroups, the normal subgroups of the full
symmetric group $S(\mu)$ on an infinite cardinal $\mu$ are all of the form $S_\kappa(\mu) =$ the
subgroup consisting of elements whose support has cardinality $< \kappa$, for some $\kappa \le
\mu^+$. A many-sorted structure ${\cal M}_{\kappa \lambda \mu}$ is defined which, it is shown,
encapsulates the first order properties of the group $S_\lambda(\mu)/S_\kappa(\mu)$. Specifically,
these two structures are (uniformly) bi-interpretable, where the interpretation of  ${\cal M}_{\kappa
\lambda \mu}$ in $S_\lambda(\mu)/S_\kappa(\mu)$ is in the usual sense, but in the other direction is
in a weaker sense, which is nevertheless sufficient to transfer elementary equivalence. By
considering separately the cases $cf(\kappa) > 2^{\aleph_0}, cf(\kappa) \le 2^{\aleph_0} < \kappa,
\aleph_0 < \kappa < 2^{\aleph_0}$, and $\kappa = \aleph_0$, we make a further analysis of the first
order theory of $S_\lambda(\mu)/S_\kappa(\mu)$, introducing many-sorted second order structures
${\cal N}^2_{\kappa \lambda \mu}$, all of whose sorts have cardinality at most $2^{\aleph_0}$, and in
terms of which we can completely characterize the elementary theory of the groups
$S_\lambda(\mu)/S_\kappa(\mu)$.
\end{abstract}

\section{Introduction}
In \cite{McKenzie}, \cite{Pinus},  \cite{Shelah1}, 
and \cite{Shelah2} a study was made of the elementary theory of
infinite symmetric groups, and a  number of natural questions arising were answered. In this paper we
examine the quotients of normal subgroups of infinite symmetric groups in the light of similar
questions. Now the normal subgroups of infinite symmetric groups are easily describable in terms of
the cardinalities of support sets. More exactly, the {\it support} of $g \in {\rm Sym}(\Omega$) is the
set of elements of $\Omega$ moved by $g$. The non-trivial normal subgroups of Sym($\mu$) where $\mu$
is an infinite cardinal  are then of the form $S_\kappa(\mu) = \{g \in {\rm Sym}(\mu): |supp
\hspace{.1cm} g| < \kappa \}$ for some cardinal $\kappa$, and the alternating group $A(\mu)$,  (see
\cite{Scott} for example), and the objects of study here are the factors
$S_\lambda(\mu)/S_\kappa(\mu)$ for $\kappa < \lambda$.

The problem of which of these groups are isomorphic is mentioned in \cite{Scott}, but we concentrate
exclusively here on the situation with regard to elementary equivalence. We shall find that many of 
the ideas from \cite{Shelah1} carry through, though with more complicated proofs.

Interpretability results about the groups $S_\lambda(\mu)/S_\kappa(\mu)$ also provide information 
about their outer automorphisms, as was explored for instance for ${\rm
Sym}(\omega)/S_\omega(\omega)$ in \cite{Alperin}, and in a related context in \cite{Giraudet}. A
survey of this aspect is given in \cite{Truss}. The result proved in \cite{Alperin} is that the
outer automorphism group of ${\rm Sym}(\omega)/ S_\omega(\omega)$ is infinite cyclic, with a  typical
outer automorphism being induced by the map $n \mapsto n+1$. The method used there incorporates a
second order interpretation of the relevant ring of sets in the quotient group. One of things we are
able to show here is that this can actually be done in a first order fashion. The existence of this
interpretation is also applied in \cite{Truss} to show that the outer automorphism group of ${\rm
Sym}(\mu)/S_\omega(\mu)$ is infinite cyclic for any $\mu$, extending the result from \cite{Alperin}.
What the outer automorphism group of $S_\lambda(\mu)/S_\kappa(\mu)$ is in general is still open---it
seems conceivable that it is trivial whenever $\kappa > \aleph_0$.  

The first order interpretation of the ring of sets in the quotient group was carried out originally by
Rubin in \cite{Rubin1} (see also \cite{Rubin2}) by a different method. Two of his main results are
\cite{Rubin1} Theorems 4.2, 4.3 which state: 

$$(1)\qquad\qquad \{(S_\lambda(\mu)/S_\kappa(\mu), {\bf B}_{\kappa \lambda \mu}; \ldots): \kappa < \lambda \le
\mu^+\}$$
is interpretable in
$$(2)\qquad\qquad \{(S_\lambda(\mu)/S_\kappa(\mu), \circ): \kappa < \lambda \le
\mu^+\}$$
where if ${\cal P}_\lambda(\mu), {\cal P}_\kappa(\mu)$ are the rings of subsets of $\mu$ of
cardinality $< \lambda, < \kappa$ respectively, then ${\bf B}_{\kappa \lambda \mu}$ is the boolean
algebra generated by ${\cal P}_\lambda(\mu)/{\cal P}_\kappa(\mu)$ and $\ldots$ signifies the `natural
relations and functions', 

(2)$$(2)\{(S(\mu)/S_\kappa(\mu): \kappa \le \mu, cf(\kappa) > 2^{\aleph_0}\}$$
is bi-interpretable (in a suitable sense) with 
$$\{\{(2^{\aleph_0},[\kappa,\mu],{\cal P}_{(2^{\aleph_0})^+}((2^{\aleph_0} \cup [\kappa,\mu])^2);
<,E)\}: \kappa \le \mu, cf(\kappa) > 2^{\aleph_0}\}$$
where $[\kappa,\mu] = \{\nu: \nu$ is a cardinal and $\kappa \le \nu \le \mu\}$, and $(a,b,R) \in E
\Leftrightarrow (a,b) \in R \in {\cal P}_{(2^{\aleph_0})^+}((2^{\aleph_0} \cup [\kappa,\mu])^2)$. 

Our
corresponding results are Theorems \ref{2.6}, \ref{4.3}, and Corollaries \ref{3.9},
\ref{4.4}.  Combining Rubin's
result with (1) and (2) with \cite{Shelah1},\cite{Shelah2}gives a full classification of the elementary types of the groups in the class 
$\{(S(\mu)/S_\kappa(\mu): \kappa \le \mu, cf(\kappa) > 2^{\aleph_0}\}$
, that is two cardinls.
The case $cf(\kappa) \le 2^{\aleph_0}$ which Rubin gives as an open 
question is treated in our final
section.

The notation used is fairly standard. We use $\kappa, \lambda$, $\mu$, and $\nu$ to stand for 
cardinals (usually infinite), and $|X|$ for the cardinality of the set $X$. If $\Omega$ is any set 
we write Sym($\Omega$) for the group of all permutations of $\Omega$ (1--1 maps from $\Omega$ 
onto itself), with permutations acting on the right, and we write $S(\mu)$ for Sym($\mu$) for
any cardinal $\mu$. For $g \in {\rm Sym} (\Omega)$  we let $supp \hspace{.1cm} g$ be the support of
$g$. If we are working in $S_\lambda(\mu)/S_\kappa(\mu)$ (where $S_\lambda(\mu),S_\kappa(\mu)$ are
as introduced above) then we refer to sets of cardinality less than $\kappa$ as {\it small}. We use
overlines such as $\overline x$ to stand for finite sequences (`tuples') $(x_1,x_2,\dots,x_n)$. By a
{\it  permutation representation} or {\it action} of a group $G$ we understand a homomorphism $\theta$
from $G$ into Sym($X$) for some set $X$. The representation is {\it faithful} if $\theta$ is 1--1,
it is {\it transitive} if for any $x, y \in X$ there is $g \in G$ such that $x(g \theta) = y$, and
it is {\it trivial} if its image is the trivial group.

If $X$ is a subset (or sequence of elements) of a group $G$, we let $\langle X \rangle$ 
denote the subgroup generated by $X$. If $g,h \in G$ we write $g^h$ for the conjugate $h^{-1}gh$ of 
$g$ by $h$. If $\overline g$ is a sequence of members of $G$ and $h \in G$, we write 
${\overline g}^h$ for the sequence whose $i$th entry is $g_i^h$, and if $\overline g, \overline h$ 
are sequences of members of $G$ of the same length, we let ${\overline g}\ast{\overline h}$ be the 
sequence whose $i$th entry is $g_ih_i$. If ${\overline g}_1^h = {\overline g}_2$ for some $h, 
{\overline g}_1$ and ${\overline g}_2$ are said to be {\it conjugate}. If $N \le G$ and $\overline f 
= (f_1,\dots,f_n) \in G^n$ we let $N.\overline f = (Nf_1,\dots,Nf_n)$.

We write ${\cal P}(X)$ for the power set of the set $X$, and ${\cal P}_\kappa(X)$ for the set of
subsets of $X$ of cardinality less than $\kappa$. Then  ${\cal P}(X)$ is a boolean algebra,
and each ${\cal P}_\kappa(X)$ for $\kappa$ infinite is a ring of sets. Moreover, if $\aleph_0
\le \kappa < \lambda \le |X|^+, {\cal P}_\kappa(X)$ is an ideal of ${\cal P}_\lambda(X)$, so we may
study the quotient ring  ${\cal P}_\lambda(X)/{\cal P}_\kappa(X)$, which is a boolean algebra just in
the case where $\lambda = |X|^+$ (that is, where ${\cal P}_\lambda(X) = {\cal P}(X)$).

\vspace{.1in}

In the remainder of this introductory section we give an outline of the main arguments of the paper.

Our analysis of the quotient groups $S_\lambda(\mu)/S_\kappa(\mu)$ is carried out using certain many
sorted structures ${\cal M}_{\kappa \lambda \mu}$ and ${\cal N}^2_{\kappa \lambda \mu}$. (There is
also a simpler version ${\cal M}^*_{\kappa \lambda \mu}$ of ${\cal M}_{\kappa \lambda \mu}$
applicable just in the case $cf(\kappa) > 2^{\aleph_0}$.) These structures are devised with the
object of describing the permutation action of tuples of elements of $S_\lambda(\mu)$, modulo
small sets. The essential properties of such an $n$-tuple ${\overline g} = (g_1, g_2, \ldots, g_n)$
are described by its action on the orbits of the subgroup $\langle {\overline g} \rangle$. In fact, if
${\overline g}_1$ and ${\overline g}_2$ are $n$-tuples of elements of $S_\lambda(\mu)$ then
${\overline g}_1$ and ${\overline g}_2$ are conjugate if and only if the  orbits of $\langle
{\overline g}_1 \rangle$ and $\langle {\overline g}_2 \rangle$ can be put into 1--1 correspondence in
such a way that the action of ${\overline g}_1$ on each orbit of $\langle {\overline g}_1 \rangle$ is
isomorphic to that of ${\overline g}_2$ on the corresponding orbit of $\langle {\overline g}_2
\rangle$. Similar remarks apply in the quotient group, except that we have to allow fewer than
$\kappa$ `mistakes' (by passing to equivalence classes of a suitable equivalence relation).

These considerations lead us to observe that what should represent ${\overline g}$ in ${\cal
M}_{\kappa \lambda \mu}$ is a list of how many $\langle {\overline g} \rangle$-orbits there are of
the various possible isomorphism types, where by `isomorphic' here we mean `under the action of
${\overline g}$'. Included among the sorts of ${\cal M}_{\kappa \lambda \mu}$ are therefore, for
each positive integer $n$, the family $IS_n$ of isomorphism types of pairs $(A, {\overline f})$,
where $\overline f$ is an $n$-tuple of permutations of $A$ acting transitively on $A$. We keep
track of the `list' of how many orbits there are of the various types by means of a function $h$
from $IS_n$ to cardinals, and the family of all these forms a further collection of sorts $F_n$. In
$F_n$ we have to identify two functions under an equivalence relation ${\cal E}_n$ if they arise from
members of $S_\lambda(\mu)$ lying in the same coset of $S_\kappa(\mu)$.

Already it is clear that ${\cal M}_{\kappa \lambda \mu}$ will have second order features (not
surprisingly, since elements of $S_\lambda(\mu)$ are {\em subsets} of $\mu^2$), but it is still
officially construed at this stage as a first order structure. The main reason for this is that at
present we cannot identify the elements of $F_n$ as functions from $IS_n$ to $Card$ (= the set of
cardinals $< \lambda$), as we would like, because, as just remarked, the members of $F_n$ are
${\cal E}_n$-classes, and ${\cal E}_n$ is not in general compatible with application. This point
is responsible for many of the complications in the paper. In the special case $cf(\kappa) >
2^{\aleph_0}$, we {\em can} so identify them, and the analysis is considerably simplified. If we
do not assume $cf(\kappa) > 2^{\aleph_0}$, then the best we can do to point the connection
between $IS_n$ and $F_n$ is to consider an `application' function $App_n$ which acts on $F_n \times
IS_n$ and gives values in $Card^- = \{ \nu \in Card: \nu = 0 \vee \kappa \le \nu < \lambda\}$. This
then {\em will} be compatible with ${\cal E}_n$, which is why all values $< \kappa$ are replaced by
0. The final sort in ${\cal M}_{\kappa \lambda \mu}$ is therefore $Card^-$, and various relations
and functions are included in its signature to express which of its properties mirror the first
order properties of $S_\lambda(\mu)/S_\kappa(\mu)$. The most important of these is $App_n$, but we
also need relations $Eq$ and $Prod$ corresponding to `equality' and `product in the group', and 
`projections' $Proj_n$ to handle existential quantification. Here $Eq \subseteq F_2$, $Prod
\subseteq F_3$, and $Proj_n$ is a function from $F_{n+1}$ to $F_n$. Corresponding relations
$Eq^1$, $Prod^1$, and $Proj_n^1$ are defined on the $IS_n$, which in `nice' cases are sufficient
to express $Eq$, $Prod$, and $Proj_n$.

The minimum goal in defining the structures ${\cal M}_{\kappa \lambda \mu}$ is that
$S_{\lambda_1}({\mu_1})/S_{\kappa_1}({\mu_1})$ and $S_{\lambda_2}({\mu_2})/S_{\kappa_2}({\mu_2})$
should be elementarily equivalent if and only if ${\cal M}_{\kappa_1 \lambda_1 \mu_1}$ and ${\cal
M}_{\kappa_2 \lambda_2 \mu_2}$ are (Corollary \ref{4.4}), and in a sense this `solves the
problem' of which of the quotient groups are elementarily equivalent. More precise information
is however avalaible. In particular, ${\cal M}_{\kappa \lambda \mu}$ is `explicitly
interpretable' in $G = S_\lambda(\mu)/S_\kappa(\mu)$; this is `interpretability' in the usual
sense, meaning that each sort and relation and function of ${\cal M}_{\kappa \lambda \mu}$ can be
represented by a definable (without parameters, in fact) relation on some power of $G$. In the
other direction we cannot hope for explicit interpretability, as one sees just by looking at the
cardinalities of the structures; a weaker property which we call `semi-interpretability'
(Definition \ref{2.5}) is established here, which is still strong enough to transfer elementary
equivalence. The fact that $S_\lambda(\mu)/S_\kappa(\mu)$ is semi-interpretable in ${\cal M}_{\kappa
\lambda \mu}$ is shown in Theorem \ref{2.6}, and essentially involves making precise the discussion
in the previous paragraph. It goes by induction on formulae of the language of group theory. For the
basis cases we use $Eq$ and $Prod$, and for the key induction step (existential quantification),
$Proj_n$.

The method for interpreting ${\cal M}_{\kappa \lambda \mu}$ in $G = S_\lambda(\mu)/S_\kappa(\mu)$ is
described in Sections 3 and 4. In Section 3 we show how the quotient ring of sets ${\cal
P}_\lambda(\mu)/{\cal P}_\kappa(\mu)$ can be interpreted. The ideas behind McKenzie's
corresponding calculations for the symmetric group \cite{McKenzie} are followed, but with considerably
greater complications. The key point is to express disjointness  of supports (`almost disjointness'
actually, meaning that they intersect in a small set). Now clearly, if two permutations have almost
disjoint supports, then they commute in $G$. The converse is very far from true, but we follow this
as a first idea, and study the configurations of certain commuting elements in sufficient details to
express disjointness. Specifically we consider sequences $\overline g$ of length 60 which satisfy the
diagram (which we write $alt_5$) of $A(5)$, the alternating group on 5 symbols, in some fixed
enumeration. This group is chosen because it is simple, and its outer automorphisms and transitive
permutation representations are easy to describe. Now apart from a small set, any 60-tuple satisfying
$alt_5$ is determined up to conjugacy by how many orbits it has of the (finitely many) possible
transitive permutation representations. Indeed this is precisely the information given by the element
of $F_{60}$ corresponding to such a tuple. By means of a (rather technical) analysis of how these
interact we can derive a formula which holds for two elements satisfying $alt_5$ if and only if (they
have a special form and) their supports are almost disjoint. Using this we find another formula which
says that two involutions have almost disjoint supports, and elements of ${\cal P}_\lambda(\mu)/{\cal
P}_\kappa(\mu)$ are then represented by (cosets of) involutions of $G$. This gives the interpretation
of the quotient ring of sets in $G$, and that of the action of $G$ on ${\cal P}_\lambda(\mu)/{\cal
P}_\kappa(\mu)$ follows easily. Moreover, all the other items of the signature of ${\cal M}_{\kappa
\lambda \mu}$ can be interpreted without much further difficulty, though there are some slight
complications in special cases, such as $\lambda = \mu^+$ or $\kappa = \aleph_0$. It is important that
we can distinguish each special case by a first order formula. For instance, the structures in which
$\lambda = \mu^+$ may be singled out by a formula saying that there is a group element such that the
only element disjoint from it is the identity (an element which moves every element of $\mu$ for
instance).

Although we generally expect ${\cal M}_{\kappa \lambda \mu}$ to have much smaller cardinality
than $S_\lambda(\mu)/S_\kappa(\mu)$, and it expresses the structure of the group in a more
compact form, the ${\cal M}_{\kappa \lambda \mu}$ still form a proper class, in view of the
presence of the sort $Card^-$. In Sections 5 and 6 we introduce the structures ${\cal N}^2_{\kappa
\lambda \mu}$, all of whose sorts have cardinality $\le 2^{\aleph_0}$, in order to be able to
reduce the problem about elementary equivalence of the groups to questions about ordinals of
cardinality $\le 2^{\aleph_0}$. In addition, the fact that ${\cal M}_{\kappa \lambda \mu}$ is a
`second order structure in disguise' is brought more out into the open, since ${\cal N}^2_{\kappa
\lambda \mu}$ genuinely {\em is} second order (hence the superscript 2). The language used to
describe ${\cal N}^2_{\kappa \lambda \mu}$ has first order variables ranging over each of its sorts,
and for each $n$-tuple of sorts, $n$-ary relations whose $i$th place lies in the $i$th sort in
the list. See Definitions \ref{5.4} and \ref{6.2}. In some cases we have to restrict the
cardinality of the relations over which the second order variables range.

Looking first at the more straightforward case, to indicate the main ideas, suppose that $cf
(\kappa) > 2^{\aleph_0}$. We show that now $App_n$ can genuinely be construed as `application', so
that we may fully describe $F_n$ in terms of $IS_n$ and $Card^-$. What therefore controls the
structure ${\cal M}_{\kappa \lambda \mu}$ is $Card^-$, and more specifically its order-type
$\alpha = \alpha(\kappa,\lambda,\mu)$. The crucial ordinals needed to describe the elementary
theory of ${\cal M}_{\kappa \lambda \mu}$ are found by writing $\alpha$ in `base $\Omega$ Canotor
normal form' where $\Omega = (2^{\aleph_0})^+$, and the countable list of ordinals $\alpha_{[n]}$
(the Cantor coefficients) and certain cofinalities $\alpha^{[n]}$ are what replace $Card^-$ in   
${\cal N}^2_{\kappa \lambda \mu}$. Theorem \ref{5.5} asserts that ${\cal N}^2_{\kappa \lambda
\mu}$ is (explicitly) interpretable in a reduct ${\cal M}^*_{\kappa \lambda \mu}$ of ${\cal
M}_{\kappa \lambda \mu}$ (and hence in ${\cal M}_{\kappa \lambda \mu}$). Non-empty subsets of
$Card^-$ of cardinality $< \Omega$ may be encoded by members of $F_2$ (we have to use $F_2$
rather than $F_1$ since $|IS_1| = \aleph_0$ but $|IS_2| = 2^{\aleph_0}$), and it is not hard to
express all the individual terms of the base $\Omega$ Cantor normal form for $\alpha$. To
express facts about cofinalities we have to quantify over binary relations on $Card^-$ of
cardinality $< \Omega$, which may be encoded using members of $F_2^2$. To express the full second
order logic described above we use longer tuples from possibly higher $F_n$s.

The transfer of properties from ${\cal N}^2_{\kappa \lambda \mu}$ to ${\cal M}_{\kappa \lambda
\mu}$ (${\cal M}^*_{\kappa \lambda \mu}$ actually suffices in this case) is not even by a
semi-interpretation. Theorem \ref{5.9} shows directly how to express an arbitrary formula
of the first order language of ${\cal M}^*_{\kappa \lambda \mu}$ by a second order formula of
the language of ${\cal N}^2_{\kappa \lambda \mu}$. Parameters are transferred using
`$k$-representations', where this means that a tuple of elements of ${\cal M}^*_{\kappa \lambda
\mu}$ (of possibly varying sorts) is represented by a (longer) tuple of elements of ${\cal
N}^2_{\kappa \lambda \mu}$ including partial maps from $IS_2$ to $IS_2$ encoding
$\alpha_{[0]}, \ldots, \alpha_{[k-1]}$ and $\alpha^{[0]}, \ldots, \alpha^{[k-1]}$.

If $cf(\kappa) \le 2^{\aleph_0}$, we can additionally interpret in ${\cal M}_{\kappa \lambda \mu}$ 
the base $\Omega$ Cantor normal form coefficients and cofinalities of the least ordinal $\alpha^*$
such that $(\exists \gamma)(\beta = \gamma + \alpha^*)$ where $\kappa = \aleph_\beta$, so this
information needs to be added to ${\cal N}^2_{\kappa \lambda \mu}$, which now includes
$\alpha^*_{[0]}, \ldots, \alpha^*_{[k-1]}$ and $\alpha^{*[0]}, \ldots, \alpha^{*[k-1]}$ (and also
$cf(\kappa)$) as additional sorts. As remarked above, since we do not now automatically know that
$\lambda > 2^{\aleph_0}$, we have to restrict the second order variables of ${\cal N}^2_{\kappa
\lambda \mu}$ to range over relations of cardinality $< \lambda$. There are some additional
complications in the cases $\kappa \le 2^{\aleph_0}$ and $\kappa = \aleph_0$, though in all cases the
outline described in the previous two paragraphs provides the basis of our analysis. Since the
precise definition of ${\cal N}^2_{\kappa \lambda \mu}$ depends on which of these cases applies, it
is important that they can all be distinguished by elementary formulae.

In summary the main conclusions are as follows. There are first order formulae of the language of
group theory distinguishing those $S_\lambda(\mu)/S_\kappa(\mu)$ for which $\lambda \le \mu$ or
$\lambda = \mu^+$, and also the cases $cf(\kappa) > 2^{\aleph_0}$, $cf(\kappa) \le 2^{\aleph_0} <
\kappa$, $\aleph_0 < \kappa \le 2^{\aleph_0}$, and $\kappa = \aleph_0$. In the case $\lambda \le
\mu$ and $cf(\kappa) > 2^{\aleph_0}$ the following holds:

for any given ordinals $\alpha_l, \alpha^l < \Omega$ there is a first order theory $T$ in the
language of group theory such that
$$  \begin{array}{c}
\mbox{if $\kappa = \aleph_\beta, \lambda = \aleph_\gamma$, $\beta + \alpha = \gamma$, and
$\alpha_{[n]} = \alpha_n, \alpha^{[n]} = \alpha^n$ for each $n$, } \\
\mbox{then the first order theory of the group $S_\lambda(\mu)/S_\kappa(\mu)$ is equal
to $T$,}
     \end{array}  $$
with similar statements in the other cases (including reference to the $\alpha^*_{[n]}, \alpha^{*[n]}$
and so on corresponding to the exact definition of ${\cal N}^2_{\kappa \lambda \mu}$).      

\section{The Basic Machinery}

Since we are aiming at a two-way interpretation, where the technically most involved step is the
representation of many notions inside the quotient group $S_\lambda(\mu)/S_\kappa(\mu)$, we
describe in this section the structure whose bi-interpretability with this group is to be shown. In
one  direction this is interpretability in the usual sense (called `explicit interpretability' in
\cite{Shelah1}), but in the other only what we may term `semi-interpretability',---which is still
sufficient for the transfer of elementary properties. We suppose that $\aleph_0 \le \kappa < \lambda
\le \mu^+$. The interpretation is most straightforward when $cf(\kappa) > 2^{\aleph_0}$, but we can
handle the general case at the expense of some additional work. In the main presentation we assume
$\kappa > \aleph_0$, indicate how the argument simplifies when $cf(\kappa) > 2^{\aleph_0}$, and what
extra is required when $\kappa = \aleph_0$. We remark that in \cite{Rubin1} Rubin showed how to
interpret the quotient ring ${\cal P}_\lambda(\mu)/{\cal P}_\kappa(\mu)$ in the group, which is also
one of our main goals, though his methods were very different from those we use. 

\begin{definition} (i) For a finite sequence $\overline f = (f_1,f_2,\dots,f_n)$ of members of
$S_\lambda(\mu)$ we let $supp \hspace{.1cm} \overline f = \bigcup_{i=1}^n supp \hspace{.1cm} f_i$.

(ii) For a positive integer $n$ let $IS_n$ be the family of isomorphism classes of pairs 
$(A,\overline g)$ where $\overline g \in ({\rm Sym}(A))^n$ and $\langle \overline g \rangle$ acts
transitively on $A$ (and if $\lambda \le \mu$, then not every $g_i$ is equal to the identity).

(iii) $Card$ = $\{ \nu: \nu$ a cardinal such that $\nu < \lambda \}$.

(iv) $Card^- = \{ 0 \} \cup \{ \nu \in Card: \kappa \le \nu \}$.

(v) If $\overline f \in (S_\lambda(\mu))^n$ let $\chi = Ch_{\overline f}$ be the function from $IS_n$
to $Card$ given by $Ch_{\overline f}((A,{\overline g})_{\cong}) = $ the number of orbits $B$ of
$\langle \overline f \rangle$ such that $(B,\overline f) \cong (A,\overline g)$.

(vi) $Ch_n = \{ Ch_{\overline f}: \overline f \in (S_\lambda(\mu))^n \}$. 

(vii) For cardinals $\kappa_1 \le \kappa_2$ we define $\kappa_2 - \kappa_1$ to be the least
cardinal $\kappa_3$ such that $\kappa_1 + \kappa_3 = \kappa_2$, and we let $|\kappa_1 - \kappa_2| = 
|\kappa_2 - \kappa_1| = \kappa_2 - \kappa_1$.

(viii) We define an equivalence relation ${\cal E}_n$ on $Ch_n$ by letting $\chi_1 {\cal E}_n
\chi_2$ if

\noindent $\sum \{ |\chi_1(t)-\chi_2(t)|: t \in IS_n \} < \kappa$. 

(ix) For each $n \ge 1$ let $F_n$ be the set of functions $h: IS_n \rightarrow Card$ such that 

$$\sum \{ h(t): t \in IS_n \} < \lambda, \mbox{ and if } \lambda = \mu^+, \sum \{ h(t): t \in
IS_n \} = \mu,$$ 

\noindent $modulo \: {\cal E}_n$.          \label{2.1} \end{definition}

\begin{remark} In `nice' cases $cf(\kappa) > 2^{\aleph_0}$ we may replace $Card$ by $Card^-$, and
then the definition of $Ch_{\overline f}$ is modified by replacing all values less than $\kappa$ by
$0$. Corresponding to this, for $\chi \in Ch_n$ we let $\chi^-(t) = \chi(t)$ if $\chi(t) \ge \kappa$
and $\chi^-(t) = 0$ otherwise. We consider this case further in section 5. If $\kappa = \aleph_0$, 
the definition of ${\cal E}_n$ is modified; here we let $\chi_1 {\cal E}_n \chi_2$ if $\sum \{
|A_t|.|\chi_1(t)-\chi_2(t)|: t \in IS_n \} < \aleph_0$ where $t = (A_t,{\overline g})_{\cong}$ (and
this says that $\chi_1(t) = \chi_2(t)$ whenever $A_t$ is infinite, and $\{ t: \chi_1(t) \neq
\chi_2(t)\}$ is finite).               \label{2.2}    \end{remark}
 
\begin{lemma} (i) For any $\overline f \in (S_\lambda(\mu))^n, \sum \{ Ch_{\overline f}(t): t \in
IS_n \} < \lambda$.

(ii) For any function $h$ from $IS_n$ to $Card$ such that $\sum \{ h(t): t \in IS_n \} < \lambda$,
and if $\lambda = \mu^+, \sum \{ h(t): t \in IS_n \} = \mu$, there is $\overline f \in
(S_\lambda(\mu))^n$ with $Ch_{\overline f} = h$. 

(iii) For ${\overline f}, {\overline g} \in (S_\lambda(\mu))^n, Ch_{\overline f} \: {\cal E}_n \:
Ch_{\overline g}$ if and only if $S_\kappa(\mu).{\overline f}$ and $S_\kappa(\mu).{\overline g}$ are
conjugate.     \label{2.3}  \end{lemma}

\noindent{\bf Proof}\hspace{.1in} (i) For each $t \in IS_n, Ch_{\overline f}(t)$ is the number of orbits of
$\langle \overline f \rangle$ on $\mu$ of that isomorphism type. Hence if $\lambda \le \mu,
\sum_tCh_{\overline f}(t)$ is equal to the number of non-trivial orbits of $\langle \overline f
\rangle$ on $\mu$, which has cardinality at most $|supp \; \overline f|$. But $supp \; \overline f$
has cardinality less than $\lambda$. If $\lambda = \mu^+,\sum_tCh_{\overline f}(t)$ equals the number
of orbits of $\langle \overline f \rangle$ on $\mu$, which is $\le \mu < \lambda$. (If $\kappa =
\aleph_0$, instead we have $\sum \{|A_t|.Ch_{\overline f}(t): t \in IS_n \} < \lambda$.)

(ii) First suppose $\lambda \le \mu$. For each $t \in IS_n$ choose a representative $(A_t,\overline
g_t)$ of that isomorphism type. Identify $\bigcup \{ A_t \times h(t): t \in IS_n \}$ with a subset
$A$ of $\mu$, and let $\overline f$ act on $A_t \times h(t)$ as $\overline g_t$ does and fix all
points outside $A$. Then for each $t, \langle \overline f \rangle$ has precisely $h(t)$ orbits of
type $t$. If $\lambda = \mu^+$ and $\sum \{h(t): t \in IS_n \} = \mu$, where now the trivial
isomorphism type is allowed, we may identify $\bigcup \{ A_t \times h(t): t \in IS_n \}$ with the
whole of $\mu$. (If $\kappa = \aleph_0, h(t)$ is replaced by $|A_t|.h(t)$.)  

(iii) Altering a member of $(S_\lambda(\mu))^n$ on a set of cardinality $< \kappa$ does not change
its $Ch_n-$value ($modulo \: {\cal E}_n$), so if $S_\kappa(\mu).{\overline f}^h = 
S_\kappa(\mu).{\overline g}$ it follows that $Ch_{{\overline f}^h} \: {\cal E}_n \: Ch_{\overline g}$.
But $h$ furnishes an isomorphism of $(A,{\overline f})$ to $(Ah,{\overline f}^h)$ for each orbit
$A$ of $\langle {\overline f} \rangle$, and so $Ch_{\overline f} = Ch_{{\overline f}^h} \: {\cal E}_n
\: Ch_{\overline g}.$

Conversely, if $Ch_{\overline f} \: {\cal E}_n \: Ch_{\overline g}$, by altering $\overline f$ on a
set of cardinality less that $\kappa$ we may suppose that $Ch_{\overline f} = Ch_{\overline g}$.
There is therefore a 1--1 correspondence between the orbits of $\langle {\overline f} \rangle$ and
$\langle {\overline g} \rangle$ which preserves the isomorphism type in $IS_n$, and which maps
singleton orbits to singleton orbits. Moreover this may be chosen having support of size $< \lambda$
(since $|supp \: {\overline f}|, |supp \: {\overline g}| < \lambda$). This gives rise to the desired
conjugacy $h$.  $\Box$

\vspace{.1in}

\begin{definition} Given the cardinals $\kappa, \lambda, \mu$ we form a many-sorted structure ${\cal
M} = {\cal M}_{\kappa \lambda \mu}$ with sorts grouped as follows:

\noindent sorts $1$: a sort $IS_n$ for each $n \ge 1$   (having cardinality $2^{\aleph_0}$ for $n
\ge 2, IS_1$ of 

cardinality $\aleph_0$),

\noindent sort $2$: $Card^-$, (in which $0$ and $\kappa$, as the first two elements, are definable,

so do not need to be explicitly named),

\noindent sorts $3$: a sort $F_n$ for each $n \ge 1$.

The signature taken is as follows:

\noindent unary relations $Eq^1$ on $IS_2$ and $Prod^1$ on $IS_3$ given by 
$$  \begin{array}{c}
Eq^1 = \{ t \in IS_2: t = ((A,g_1,g_2))_{\cong} \rightarrow g_1 = g_2 \}, \\
Prod^1 = \{ t \in IS_3: t = ((A,g_1,g_2,g_3))_{\cong} \rightarrow g_1 g_2 = g_3 \},
     \end{array}  $$
for each $n$ a binary relation $Proj_n^1 \subseteq IS_{n+1} \times IS_n$ given by
$$  \begin{array}{c}
Proj_n^1 = \{ (t_1,t_2) \in IS_{n+1} \times IS_n: \exists A_1 \exists A_2 \exists g_1 \exists g_2
\ldots \exists g_{n+1} (A_1 \supseteq A_2 \\
\wedge \: t_1 = ((A_1, g_1, \ldots , g_{n+1}))_{\cong} \: \wedge \: t_2 = ((A_2, g_1
| A_2 , \ldots , g_n | A_2 ))_{\cong}) \},
     \end{array}  $$ 
$<$ on sort $2$, the usual ordering of cardinals,

\noindent for each $n$ a function $App_n$ from $F_n \times IS_n$ to $Card^-$ given by $App_n(x,y) = 
\nu$ provided that for some $h$ with $(h)_{{\cal E}_n} = x, h(y) = \nu$ (noting that the value of
$h(y)$ is well-defined for $\nu \ge \kappa$, and for $\nu < \kappa$, all values are replaced by 
$0$,---see the definition of $\chi^-$ above),

\noindent unary predicates $Eq$ on sort $F_2$ and $Prod$ on sort $F_3$ given by $Eq(h), Prod(h)$ 
hold if
$$  \begin{array}{c}
\sum \{ h(t): t \in IS_2 \: \wedge (t = ((A,g_1,g_2))_{\cong} \rightarrow g_1 \neq g_2) \} = 0, \\
\sum \{ h(t): t \in IS_3 \: \wedge (t = ((A,g_1,g_2,g_3))_{\cong} \rightarrow g_1g_2 \neq g_3)
\} = 0,
     \end{array}  $$ 
respectively, (where as the sorts $3$ consist of functions modulo ${\cal E}_n$, saying that
these sums are zero means in effect that they are $< \kappa$),

\noindent and functions $Proj_n$ from sort $F_{n+1}$ to sort $F_n$ such that if $h: IS_{n+1} 
\rightarrow Card$ then $Proj_n(h): IS_n \rightarrow Card$ is given by $Proj_n(h)(t) = \sum \{
|B_{t't}|.h(t'): t' \in IS_{n+1} \}$ where for each $t' = ((A,\overline g))_{\cong} \in IS_{n+1}$, and
$\overline g$ of length $n+1, B_{t't}$ is the set of all orbits of $\langle g_1, g_2, \ldots , g_n
\rangle$ on $A$ on which $(g_1, g_2, \ldots , g_n)$ has isomorphism type $t$. (Note that $|B_{t't}|$
is independent of the particular choice of $(A,\overline g)$ corresponding to $t'$. Note also
that strictly speaking here and in the definition of $Eq$, $Prod$, we should work with the
${\cal E}_n$-classes determined by $h, Proj_n(h)$.)      \label{2.4}     \end{definition}

We include $Proj_n$ in order to handle existential quantification in the forthcoming induction
(Theorem \ref{2.6}). The definitions of $Eq$ and $Prod$ apply just in the case $\kappa >
\aleph_0$, and are intended to express equality and products in $S_\lambda(\mu)$ up to fewer than
$\kappa$ mistakes. For $\kappa = \aleph_0$, instead of summing the relevant $h(t)$ we sum 
$|A_t'|.h(t)$ where $t = ((A_t,f_1,f_2))_{\cong}$ or $((A_t,f_1,f_2,f_3))_{\cong}$ and $A_t' =
\{ \alpha \in A_t: \alpha f_1 \neq \alpha f_2 \}$ or $\{ \alpha \in A_t: \alpha f_1 f_2 \neq \alpha
f_3 \}$ respectively.

In the general case the inclusion of the sorts $IS_n$ and $Card^-$ is unnecessary, at any rate as
far as the proof of Theorem \ref{2.6} is concerned. On the other hand in all cases they can be
naturally represented within $S_\lambda(\mu)/S_\kappa(\mu)$, and in the special case $cf(\kappa) >
2^{\aleph_0}$ the `application' functions $App_n$ genuinely identify the members of $F_n$ as
functions from $IS_n$ to $Card^-$ (since here the equivalence relation ${\cal E}_n$ can be dispensed
with), meaning that $App_n': F_n \rightarrow (Card^-)^{IS_n}$ given by $App_n'(h)(t) = App_n(h,t)$ is
1--1. In section 5 we shall also see that $Proj_n$ is definable from $Proj_n^1$ and $App_n$ in this
case, and similarly for $Eq$ and $Prod$, easing the analysis of the ${\cal M}_{\kappa \lambda \mu}$.

The sense in which we can show that $S_\lambda(\mu)/S_\kappa(\mu)$ is interpretable in ${\cal
M}_{\kappa \lambda \mu}$ is weaker than the usual one and is given in the following definition.

\begin{definition} For structures $\cal M$ and $\cal N$ we say that $\cal M$ is {\em
semi-interpretable} in $\cal N$ if there is a recursive function $F$ from formulae of the language
of $\cal M$ to formulae of the language of $\cal N$ and there are functions $f_n: {\cal M}^n
\rightarrow {\cal N}$ such that for all ${\overline a} \in {\cal M}^n$ and $\varphi({\overline
x})$ with $n$ free variables, ${\cal M} \models \varphi[{\overline a}] \Leftrightarrow {\cal N}
\models F(\varphi)[f_n({\overline a})]$. If the same $F$ serves over a class of pairs of structures
then we say that the first of each pair is {\em uniformly semi-interpretable} in the second.
\label{2.5} \end{definition}

\begin{theorem} For every first order formula $\varphi(x_0, \ldots ,x_{n-1})$ of the theory of
groups there is an effectively determined first order formula $\psi(y)$ of the language of $\cal M$
such that for all $\kappa, \lambda, \mu$, and for every ${\overline f} \in (S_\lambda(\mu))^n$,

$$S_\lambda(\mu)/S_\kappa(\mu) \models \varphi[S_\kappa(\mu).{\overline f}] \Leftrightarrow {\cal
M}_{\kappa \lambda \mu} \models \psi[(Ch_{\overline f})_{{\cal E}_n}].$$   \label{2.6} \end{theorem}

\noindent{\bf Proof}\hspace{.1in} We construct $\psi$ by induction. First suppose that $\varphi$ is atomic.
It suffices to consider formulae of the form $x_0 = x_1$ and $x_0 x_1 = x_2$ for variables $x_0,
x_1, x_2$. If $\varphi(x_0,x_1)$ is $x_0 = x_1$ we take for $\psi(y)$ the formula $Eq(y)$. Then

\vspace{.1in}
\noindent $S_\lambda(\mu)/S_\kappa(\mu) \models S_\kappa(\mu).f_1 = S_\kappa(\mu).f_2$
\vspace{-.1in}
\begin{eqnarray*} & \Leftrightarrow & | \{ \alpha:\alpha f_1 \neq \alpha f_2 \} | < \kappa \\  
& \Leftrightarrow & \mbox{ the union of the orbits of } \langle f_1,f_2 \rangle \mbox{ on which the
       actions of } f_1 \mbox{ and } f_2 \\
& & \mbox{ are distinct has cardinality } < \kappa \hspace{.3in} \mbox{(since } \kappa > \aleph_0) \\
& \Leftrightarrow & \sum \{ |A_t|.Ch_{\overline f}(t): t \in IS_2 \: \wedge (t =
       ((A_t,g_1,g_2))_{\cong} \rightarrow g_1 \neq g_2) \} < \kappa  \hspace{1in} \\ 
& \Leftrightarrow & {\cal M} \models Eq((Ch_{\overline f})_{{\cal E}_n})  \hspace{1in} \mbox{(using }
\kappa > \aleph_0 \mbox{ again)} \\ 
& \Leftrightarrow & {\cal M} \models \psi[(Ch_{\overline f})_{{\cal E}_n}],
\end{eqnarray*}
and similarly for the formula $x_0x_1 = x_2$ (using $Prod$).

The propositional induction steps are straightforward. 

Finally suppose that $\varphi(\overline x)$ is $\exists y \varphi_1(\overline x,y)$. 

It is easily checked that for any ${\overline f} = (f_1, f_2, \ldots , f_n)$ and $f_{n+1}$ in
$S_{\lambda}(\mu),$ 
$$Proj_n(Ch_{{\overline f},f_{n+1}}) = Ch_{\overline f}.$$ 
>From this it follows that if $h \in F_{n+1}$, then $Proj_n(h) = Ch_{\overline f}$ if and only if $h =
Ch_{{\overline f},f_{n+1}}$ for some $f_{n+1}$. Continuing the proof we deduce that
 
\vspace{.1in}
\noindent $S_\lambda(\mu)/S_\kappa(\mu) \models \exists x \varphi_1(S_\kappa(\mu).\overline f,x)$ 
\vspace{-.1in}
\begin{eqnarray*} & \Leftrightarrow & \mbox{ for some } g \in S_\lambda(\mu),
S_\lambda(\mu)/S_\kappa(\mu) \models \varphi_1[S_\kappa(\mu).\overline f,S_\kappa(\mu).g]   \\
& \Leftrightarrow & \mbox{ for some } g \in S_\lambda(\mu), {\cal M} \models 
               \psi_1[(Ch_{\overline f,g})_{{\cal E}_{n+1}}]       \\  
& \Leftrightarrow & {\cal M} \models \exists x (\psi_1(x) \wedge Proj_n(x) = (Ch_{\overline
f})_{{\cal E}_n}),
\end{eqnarray*}
where $\psi_1$ is a formula corresponding to $\varphi_1$ as given by the induction hypothesis, and 
so we take for $\psi(\overline x)$ the formula $\exists y(\psi_1(y) \wedge Proj_n(y) = \overline x)$.
 $\Box$

\section{Interpreting ${\cal P}_\lambda(\mu)/{\cal P}_\kappa(\mu)$ in $S_\lambda(\mu)/S_\kappa(\mu)$}

In this section we show how it is possible to interpret many `set-theoretical' properties inside
$S_\lambda(\mu)/S_\kappa(\mu)$, by representing subsets of $\mu$ via supports of suitably chosen
elements (always up to fewer than `$\kappa$ mistakes'), and consequently to interpret the ring
${\cal P}_\lambda(\mu)/{\cal P}_\kappa(\mu)$. The key idea is to use sequences whose entries are
transitive representations of a specific finite {\it non-abelian} group to represent the subsets,
which enables us to capture disjointness of their supports via a commutativity condition. We
introduce the necessary formulae one by one, and outline why they represent what is required. 

\vspace{.1in}

Let $G$ be a fixed finite group of order $n$, and let $\overline G = (a_1,a_2,\dots,a_n)$ be a 
fixed enumeration of $G$ for which $a_1 = id$, the identity. In what follows we shall in fact just 
use $G = A(5)$, the alternating group on $\{0,1,2,3,4\}$. This is for three reasons: it is the 
smallest non-abelian simple group; its transitive permutation representations are easy to describe; 
and (a small point needed in the proof) its outer automorphism group is also well known (and is 
just $S(5)$).

Let $diag({\overline G}, \overline x$) be the conjunction over all $i,j,k$ between $1$ and $n$ for
which $a_ia_j = a_k$ of the formulae $x_ix_j = x_k$. This is intended to say that 
$(x_1,x_2,\dots,x_n)$ is a `copy' of $G$ (in the specified enumeration), but actually just says 
that it is a homomorphic image. We write $diag({\overline {A(5)}},{\overline x})$ as $alt_5(\overline
x)$.

\begin{lemma} Suppose $\overline f \in (S_\lambda(\mu))^n$ is such that $S_\lambda(\mu)/S_\kappa(\mu) 
\models diag(\overline G, S_\kappa(\mu).\overline f)$. Then there is a small union $X$ of $\langle 
\overline f \rangle$-orbits such that if $\alpha \in \mu - X$, and $a_ia_j = a_k$, then 
$\alpha f_if_j = \alpha f_k$.  \label{3.1} \end{lemma}

\noindent{\bf Proof}\hspace{.1in} Let $X_{ijk} = \{ \alpha \in \mu: \alpha f_if_j \neq \alpha f_k \}$ and let 
$X = \bigcup \{ X_{ijk}: a_ia_j = a_k \}$. By definition of $diag$, and of 
$S_\lambda(\mu)/S_\kappa(\mu), |X| < \kappa$. So it suffices to observe that $\mu - X$ is closed 
under the action of $\langle \overline f \rangle$. Let $\alpha \not \in X$ and $1 \le r \le n$.
Suppose that $i,j,k$ are arbitrary subject to $a_ia_j = a_k$. Then there are $s,t$ such that 
$a_ra_i = a_s$ and $a_ra_k = a_t$. We find that $a_sa_j = a_ra_ia_j = a_ra_k = a_t$. Since $\alpha
\not \in X, \alpha f_rf_if_j = \alpha f_sf_j = \alpha f_t = \alpha f_rf_k$. Thus $\alpha f_r \not 
\in X_{ijk}$ and so $\alpha f_r \not \in X$ as required.  $\Box$

\begin{lemma} For any $\overline f, S_\lambda(\mu)/S_\kappa(\mu) \models
alt_5(S_\kappa(\mu).{\overline f})$ if and only if there is a small union $X$ of orbits of $\langle 
\overline f \rangle$ on $\mu$ such that for every orbit $Y$ of $\langle \overline f \rangle$ on 
$\mu - X$, the action of $\langle \overline f \rangle$ on $Y$ is isomorphic to some action of
$A(5)$ (so that $|Y| = 1, 5, 6, 10, 12, 15, 20, 30$, or $60$, and then we say that $\overline f$ 
acts as $A(5)$ with this degree). \label{3.2} \end{lemma}

\noindent{\bf Proof}\hspace{.1in} This is immediate from Lemma~\ref{3.1} on remarking that for orbits $Y$ 
outside $X$, the action of $\langle \overline f \rangle$ on $Y$ precisely corresponds to some
transitive action of $A(5)$. The fact that the possible values of $|Y|$ are as stated follows from
the fact that any transitive action of $A(5)$ is isomorphic to its action on a coset space $[A(5):H]$
for some subgroup $H$ of $A(5)$, and the possible orders of subgroups of $A(5)$ are $1,2,3,4,5,6,10,
12$, and $60$.  $\Box$
 
\vspace{.1in}

With this lemma in mind we may define for any $\overline f \in S_\lambda(\mu)$ of length $60$ the
cardinals $\nu_m(\overline f)$ for $m \in \{ 1,5,6,10,12,15,20,30,60 \}$ by $\nu_m(\overline f) =$
the number of $\langle \overline f \rangle$-orbits on $\mu$ on which $\overline f$ acts as $A(5)$
with degree $m$. The significant values (that is, those which are preserved under passing to the
coset $S_\kappa(\mu).\overline f$) are those $\nu_m(\overline f)$ which are $\ge \kappa$, and
these provide a `profile' of $\overline f$ characterizing it up to conjugacy. 

To make further progress we need to analyse with some care some properties of the possible 
faithful transitive actions of $A(5)$, which we do in the next three lemmas.
 
\begin{lemma} Suppose that $H$ and $K$ are proper subgroups of $A(5)$. Then for some $a \in A(5), 
|H \cap a^{-1}Ka| \le 3$. Moreover, if there is $a$ such that $|H \cap a^{-1}Ka| = 3$ but no $b$
such that $|H \cap b^{-1}Kb| < 3$, then $|H| = |K| = 12$.   \label{3.3} \end{lemma}

\noindent{\bf Proof}\hspace{.1in} As $A(5)$ is simple, $|H|,|K| \le 12$. If $H$ or $K$ has order $\le 3$ we just
let $a = id$. Assuming without loss of generality that $|H| \ge |K|$ we are left with the following 
possibilities for $(|H|,|K|)$: 

$(12,12), (12,10), (12,6), (12,5), (12,4), (10,10), (10,6), (10,5), (10,4)$, 

\noindent $(6,6), (6,5), (6,4), (5,5), (5,4), (4,4)$.

The subgroups of $A(5)$ of orders $12,10,6,5,4$ are determined uniquely up to conjugacy in $S(5)$
(as is easy to check) and so by replacing by a conjugate by a member of $A(5)$ may be taken to lie
in the following list:

$12:  A(4)$  (regarded as the stabilizer of $4$ in $A(5)$), 

$10:  \langle(01234),(14)(23) \rangle, \langle(01243),(13)(24) \rangle$, 

$6: \hspace{.1in}  \langle(012),(01)(34) \rangle$,

$5: \hspace{.1in} \langle (01234) \rangle, \langle (01243) \rangle$, 

$4: \hspace{.1in} \langle(01)(23),(02)(13) \rangle$.

\noindent The following cases can be at once ruled out as $|H \cap K| \le 3$ is already true: 
$(12,10), (12,6), (12,5), (10,6), (10,4), (6,5), (6,4), (5,4)$. In all the remaining cases, which 
are $(12,12), (12,4), (10,10), (10,5), (6,6), (5,5)$, and $(4,4)$, the conjugator $(234)$ will 
serve as $a$, as is easy to check. 

Now for the final part, suppose that $|H \cap a^{-1}Ka| = 3$ for some $a$ and that 
$|H \cap b^{-1}Kb| \ge 3$ for all $b$. Then $|H|$ and $|K|$ are multiples of $3$. If $H$ or $K$ has 
order $12$ or $6$, we take it as above, and if $3$ we take it as $\langle(012)\rangle$. In all cases
except for $|H| = |K| = 12$ we find that for $b = (243), |H \cap b^{-1}Kb| = 1$ or $2$, and we
conclude that $H$ and $K$ must both have order $12$.     $\Box$

\begin{lemma} Let $D$ be the diagonal subgroup $\{(a_i,a_i): 1 \le i \le 60 \}$ of $A(5) \times 
A(5)$. Then for any subgroup $H$ of $A(5) \times A(5)$ of order $12$ or $36$, there is $a$ such that 
$|a^{-1}Ha \cap D| \neq 3$. \label{3.4} \end{lemma}

\noindent{\bf Proof}\hspace{.1in} Suppose otherwise. Thus $|H| = 12$ or $36$, and for every $a \in 
A(5) \times A(5), |a^{-1}Ha \cap D| = 3$. In particular $|H \cap D| = 3$ so we suppose that $H \cap D
= \langle ((012),(012)) \rangle$. Let $a = ((13)(24),id)$. Since $|a^{-1}Ha \cap D| = 3$ there are
$i < j < k$ such that $((ijk),(ijk)) \in a^{-1}Ha$. Also $((034),(012)) \in a^{-1}Ha$. If
$| \{i,j,k\} \cap \{0,3,4 \}| = 1$ then $\langle (ijk),(034) \rangle$ contains an element of order
$5$, contrary to $|a^{-1}Ha| = 12$ or $36$. Hence $| \{i,j,k\} \cap \{0,3,4 \}| = 2$ or $3$.
Similarly $|\{ i,j,k \} \cap \{ 0,1,2 \}| = 2$ or $3$. Therefore $i = 0$ and $j = 1$ or $2, k = 3$ or
$4$.

\noindent{\bf Case 1}: $(ijk) = (013)$. Then $((012),(012)),((031),(013)) \in H$. But these two
elements generate a group of order $144$ ($A(4) \times A(4)$ in fact).

\noindent{\bf Case 2}: $(ijk) = (014)$. Therefore $((012),(012)),((032),(014)) \in H$. Since $H$
has no element of order $5, H \le A(4) \times A(\{0,1,2,4\})$. Let $b = ((014),id)$. Then $b^{-1}Hb
\le A( \{ 1,4,2,3 \} ) \times A( \{ 0,1,2,4 \} )$. If $((i'j'k'),(i'j'k'))$ lies in 
$b^{-1}Hb$ with $i' < j' < k'$ then $\{ i',j',k' \} \subseteq \{1,4,2,3 \} \cap \{ 0,1,2,4 \}$, so
$(i'j'k') = (124)$. Then $((021),(124)) \in H$, so that $(id, (02)(14)) \in H$, from which it follows
that $|H| \neq 12, 36$.

\noindent{\bf Case 3}: $(ijk) = (023)$. Then $((012),(012)),((041),(023)) \in H$, and we argue as in
Case 2, with $b = ((013), id)$. This time we find that $((021),(123)) \in H$, so that $(id, (02)(13))
\in H$, and $|H| \neq 12,36$.

\noindent{\bf Case 4}: $(ijk) = (024)$. Then $((012),(012)), ((042), (024)) \in H$ so
$((01)(24),(014))$ $\in H$, and $|H| = 144$ as in Case 1.   $\Box$

\begin{lemma} Suppose that $\overline f, \overline g$ are subgroups of {\rm Sym}($X$) isomorphic
to $A(5)$ (in the specified listings) which centralize each other, and such that $\langle \overline f,
\overline g \rangle$ is transitive on $X$. Then $\overline f \ast \overline g$ has an orbit of
length at least $20$. Moreover, if $\overline f\ast \overline g$ has an orbit of length $20$ then
it also has an orbit of some other length greater than $1$. \label{3.5} \end{lemma}

\noindent{\bf Proof}\hspace{.1in} Let ${\cal X} = \{ X_i: i < m \}$ and ${\cal Y} = \{ Y_j: j < n \}$ be the 
families of orbits of $\overline f$ and $\overline g$ respectively. Then as $\overline f$ and 
$\overline g$ commute, $\overline f$ and $\overline g$ each preserve ${\cal X}$ and ${\cal Y}$
(setwise), and hence also ${\cal Z} = \{ X_i \cap Y_j: i < m, j < n \}$. Moreover by transitivity
of $\langle \overline f, \overline g \rangle$ on $X$ the actions of $\overline f$ on its orbits
are all isomorphic, as are the actions of $\overline g$ on its orbits. Since $\overline f, 
\overline g$ are isomorphic to $A(5)$, these orbits are all non-trivial, and since $A(5)$ is
simple, they all have at least $5$ members.

\noindent{\bf Case 1}: $m = n = 1$. Thus $\overline f$ and $\overline g$ are both transitive.

In this situation it is standard that $\overline f$ and $\overline g$ both act regularly (see
\cite[Theorem 3.2.9]{Tsuzuku}). For suppose that $xf_i = x$. Then for each $j, (xg_j)f_i = xf_ig_j =
xg_j$ and as $\overline g$ is transitive, $f_i = id$. Similarly $\overline g$ is regular. By suitably
labelling the elements of $X$ we may suppose that $X = A(5)$ and $\overline f$ is the right regular
action, in other words $(a_i)f_j = a_ia_j$ for each $i$ and $j$.

Now we appeal essentially to the fact that the centralizer of the right regular action is the left
regular action (see \cite[Theorem 3.2.10]{Tsuzuku}). Let $a_1g_i = a_r$. Then $a_jg_i = a_1f_jg_i =
a_1g_if_j = a_rf_j = a_ra_j$. Hence $g_i$ is multiplication on the left by $a_r$. Let us write $a_r$
as $a_i \theta$. Thus $\theta$ is 1--1 since if $a_i\theta = a_j\theta, a_1g_i = a_1g_j$ and $i=j$
(by regularity). So also $\theta$ is onto. Moreover it is an anti-homomorphism, since
$(a_ia_{i'})\theta = a_1(g_ig_{i'}) =  ((a_i\theta)a_1)g_{i'} = (a_{i'}\theta)(a_i\theta)a_1 =
(a_{i'}\theta)(a_i\theta)$. Thus $\varphi$  given by $a_i\varphi = (a_i^{-1})\theta$ is an
automorphism of $A(5)$. So for some $s \in S(5),  a_i\varphi = s^{-1}a_is$ for all $i$, so that
$a_i(f_jg_j) = s^{-1}a_j^{-1}sa_ia_j$. Now the length of the orbit of $f_jg_j$ containing $a_i$ is
equal to the index of its stabilizer in $A(5)$. But
$s^{-1}a_j^{-1}sa_ia_j = a_i \Leftrightarrow a_j^{-1}sa_i = sa_ia_j^{-1} \Leftrightarrow 
a_j \in C_{A(5)}(sa_i) = A(5) \cap C_{S(5)}(sa_i)$. Now $sa_i$ either ranges over $A(5)$ or over 
$S(5) - A(5)$. If $A(5)$ let $sa_i = (012)$ or $(01234)$. Then $|C_{A(5)}(sa_i)| = 3$ or $5$ and 
so there are orbits of lengths $20$ and $12$. If $S(5) - A(5)$ let $sa_i = (0123)$. Then 
$|C_{A(5)}(sa_i)| = 2$ and so there is an orbit of length $30$.

\noindent{\bf Case 2}: $m = 1 \hspace{.1cm} \wedge \hspace{.1cm} n > 1$ (or similarly $m > 1 
\hspace{.1cm} \wedge \hspace{.1cm} n = 1$).

Then $\overline f$ is transitive, so by the same proof as above, $\overline g$ acts semiregularly
(that is, only the identity has any fixed point). Hence $\overline g$ acts regularly on each orbit,
and so each orbit has size $60$. But then $|X| > 60$, contrary to $\overline f$ transitive on $X$.

\noindent{\bf Case 3}: $m,n > 1$.

Since $\langle \overline f, \overline g \rangle$ is transitive, the actions of $\overline f$ on $\cal
Y$ and $\overline g$ on $\cal X$ are both transitive, and hence faithful. Moreover $\langle \overline
f, \overline g \rangle$ acts transitively on ${\cal Z} = \{ X_i \cap Y_j: i < m, j < n \}$ (which in
particular means that all $X_i \cap Y_j$ are non-empty of equal size).

We show that some orbit of $\overline f\ast \overline g$ in its action on $\cal Z$ has length $\ge 
20$, and it will follow that the same applies to its action on $X$. Now the length of the orbit 
containing $X_i \cap Y_j$ is equal to the index of its stabilizer, and as $(X_i \cap Y_j)f_kg_k
= X_ig_k \cap Y_jf_k, \hspace{.1cm} (X_i \cap Y_j)f_kg_k = X_i \cap Y_j \Leftrightarrow X_ig_k = 
X_i \wedge Y_jf_k = Y_j$. Hence $\{ a_k: (X_i \cap Y_j)f_kg_k = X_i \cap Y_j \} = \{ a_k: X_ig_k = X_i \}
\cap  \{ a_k: Y_jf_k = Y_j \}$ and all we have to do is to show that for some $i,j,$ the right hand
side  has order $\le 3$. Let $H = \{ a_k: X_0g_k = X_0 \}$ and $K = \{ a_k: Y_0f_k = Y_0 \}$. Then 
the stabilizers of the other $X_i$ and $Y_j$ are just the conjugates of these. For instance 
$\{ a_k: X_0g_ig_k = X_0g_i \} = \{ a_k: a_ia_ka_i^{-1} \in H \} = a_i^{-1}Ha_i$ and 
$\{ a_k: Y_0f_jf_k = Y_0f_j \} =  a_j^{-1}Ka_j$.

For our choice we take $i = 0$ and select $j$ by using Lemma~\ref{3.3}. 

Finally we have to show (still in Case $3$) that not all orbits of $\overline f \ast \overline g$
can have length $20$ or $1$. Suppose otherwise. Since $m \ge 5$ and $\overline f \ast \overline g$
acts transitively on $\{ X_i: i < m \}$, none of the orbits can have length $1$. Applying the last 
clause of Lemma~\ref{3.3} we find that $|H| = |K| = 12$, and so $m = n = 5$. Therefore $|X| = 
25|X_0 \cap Y_0|$ and since this is a multiple of $20$, and $|X_0| = 5|X_0 \cap Y_0|$ is a factor of
$60,  |X_0 \cap Y_0| = 4$ or $12$, so $|X| = 100$ or $300$.

Pick $x \in X_0 \cap Y_0$ and let $L = \{ (a_i,a_j) \in A(5) \times A(5) : xf_ig_j = x \}$. Since
$A(5) \times A(5)$ acts transitively on $X$ via $(\overline f,\overline g), |L| = 60^2/|X| = 12$ or
$36$. By Lemma~\ref{3.4} (and with $D$ as there), there are $i,j$ such that
$|(a_i,a_j)^{-1}L(a_i,a_j)  \cap D| \neq 3$. Let $y = xf_ig_j$. Then $yf_kg_k = y \Leftrightarrow
xf_if_kg_jg_k = xf_ig_j \Leftrightarrow xf_if_kf_i^{-1}g_jg_kg_j^{-1}$

\noindent $= x \Leftrightarrow (a_ia_ka_i^{-1},a_ja_ka_j^{-1}) \in L \Leftrightarrow (a_k,a_k) \in
(a_i,a_j)^{-1}L(a_i,a_j)$. Hence $|\{a_k: yf_kg_k = y \}| = |(a_i,a_j)^{-1}L(a_i,a_j) \cap D| \neq
3$, and so the orbit of $y$ under the action of $\overline f \ast \overline g$ does not have length
$20$ after all.    $\Box$

\vspace{.1in}   

We now move towards the construction of a formula which is intended to say that $\overline x$ acts as
$A(5)$ on all but a small set of its orbits, and that each such orbit has length $1$ or $5$. Actually
we stop short of doing this (even though it can be done) and just find a formula restricting the
range of representations possible---as this provides a quicker route to our goal.  We require the
following auxiliary formulae:

\vspace{.1in}

\noindent $comm_{m,n}(\overline x, \overline y$): $\bigwedge_{1 \le i \le m,1 \le j \le n}x_iy_j =
y_jx_i$, 

\noindent \hspace{1in} where $m$ and $n$ are the lengths of $\overline x$ and $\overline y$. This 
asserts that 

\noindent \hspace{1in} each entry of $\overline x$ commutes with each entry of $\overline y$.

\noindent $conj_n(\overline x, \overline y$): $(\exists z)({\overline x}^z = {\overline y})$.

In practice we omit the subscripts from $comm_{m,n}$ and $conj_n$ (and other similar formulae).

\noindent $indec(\overline x$): $alt_5(\overline x) \wedge (\forall \overline y)(\forall \overline
z)(comm(\overline y,\overline z) \wedge alt_5(\overline y) \hspace{.1cm} \wedge$

\noindent \hspace{1in} $alt_5( \overline z) \wedge {\overline x} = {\overline y}\ast{\overline z}
\rightarrow  (conj(\overline x,\overline y) \vee conj(\overline x,\overline z))$.

\begin{lemma} For any sequence $\overline f$ of elements of $S_\lambda(\mu)$ of length $60$,
$S_\lambda(\mu)/S_\kappa(\mu)$ 

\noindent $\models indec(S_\kappa(\mu).{\overline f})$ if and only if $|\mu - supp \; {\overline f}|
= \mu$, $\langle \overline f \rangle$ acts as $A(5)$ on all orbits outside a small subset of $\mu$,
$\nu_{30}(\overline f), \nu_{60}(\overline f) < \kappa$, and there is at most one $m \in 
\{5,6,10,12,15,20 \}$ for which $\nu_m(\overline f) \ge \kappa$. \label{3.6} \end{lemma}

\noindent{\bf Proof}\hspace{.1in} We remark that we need to stipulate $|\mu - supp \; {\overline f}|
= \mu$ in view of the possibility that $\lambda = \mu^+$. Let us say that $\overline f \not \in
S_\kappa(\mu)$ is {\it indecomposable} if $S_\lambda(\mu)/S_\kappa(\mu) \models
indec(S_\kappa(\mu).\overline f)$.

First observe that if $\nu_{m_1}(\overline f), \nu_{m_2}(\overline f) \ge \kappa$ where $1 < m_1
< m_2$ then we may write $\overline f$ as $\overline g\ast\overline h$ where $\overline g$
is the restriction of $\overline f$ to the union of its orbits of length $m_1$ (that is it agrees 
with $\overline f$ there and fixes all other points), and $\overline h$ is the restriction of 
$\overline f$ to the complement of the union of these orbits. Clearly $\overline f, \overline g$
commute, $S_\kappa(\mu).\overline g$ and $S_\kappa(\mu).\overline h$ satisfy $alt_5(\overline x)$,
and $\overline f = \overline g\ast\overline h$. But neither $S_\kappa(\mu).\overline g$ nor
$S_\kappa(\mu).\overline h$ is conjugate to $S_\kappa(\mu).\overline f$.

Next suppose that $\nu_{30}(\overline f) \ge \kappa$, and let $X$ be the union of all orbits of 
$\overline f$ of length $30$ on which $\overline f$ acts as $A(5)$. We let $\overline g$ and 
$\overline f$ agree on $\mu - X$ and $\overline h$ fix $\mu - X$ pointwise. Let $Y$ be a typical
orbit of $\overline f$ contained in $X$, (and hence of length $30$). Now if $H$ and $K$ are 
subgroups of $A(5)$ of orders $12$ and $10$, then $A(5)$ has a transitive action of degree $30$ on 
$[A(5):H] \times [A(5):K]$, (since $|H \cap K| = 2$), which is therefore isomorphic to the action of
$\langle \overline f \rangle$ on $Y$. So we may let $Y = \{\alpha_{(Ha_i,Ka_j)}: Ha_i \in [A(5):H],
Ka_j \in [A(5):K] \}$ in such a way that for each $k, \alpha_{(Ha_i,Ka_j)}f_k =
\alpha_{(Ha_ia_k,Ka_ja_k)}$. The point is that this expresses  the action of $\overline f$ on $Y$ as
a commuting `product' of actions having orbits of sizes 5 and 6. We let $\alpha_{(Ha_i,Ka_j)}g_k
= \alpha_{(Ha_ia_k,Ka_j)}$ and $\alpha_{(Ha_i,Ka_j)}h_k = \alpha_{(Ha_i,Ka_ja_k)}$. This therefore
defines the actions of $\overline g$ and $\overline h$ on  the orbits of $\overline f$ having length
30. It is clear that neither $S_\kappa(\mu).\overline g$ nor $S_\kappa(\mu).\overline h$ can be
conjugate to $S_\kappa(\mu).\overline f$, since $\nu_{30}(\overline g),\nu_{30}(\overline h) <
\kappa$. But $S_\kappa(\mu).\overline g$ and $S_\kappa(\mu).\overline h$ fulfil the other requirements
on $\overline y$ and $\overline z$ in $indec$, and so we conclude that $S_\kappa(\mu).\overline f$
cannot satisfy $indec$.

If $\nu_{60}(\overline f) \ge \kappa$, a similar argument applies, but this time taking $|H| = 12$
and $|K| = 5$. 

Now suppose that $|\mu - supp \; {\overline f}| < \mu$. Let $X$ be a union of orbits of $\langle
{\overline f} \rangle$ such that $|X| = |\mu - X| = \mu$, and let $\overline g$ and $\overline h$ be
the restrictions of $\overline f$ to $X$ and $\mu - X$ respectively. Then $S_\kappa(\mu).{\overline
g}$ and $S_\kappa(\mu).{\overline h}$ provide witnesses for $y$ and $z$ violating
$indec(S_\kappa(\mu).{\overline f})$.

Conversely, suppose that $|\mu - supp \; {\overline f}| = \mu$ and for some $m \in \{5,6,10,12,
15,$ $20\}$, the union $X$ of the orbits of $\langle \overline f \rangle$ of length $m$ on which
$\overline f$ acts as $A(5)$ has cardinality $\ge \kappa$, and that $\overline f$ fixes all but a
small subset $Y$ of $\mu - X$. We verify $indec(S_\kappa(\mu).\overline f)$. Suppose
$S_\kappa(\mu).\overline g$ and $S_\kappa(\mu).\overline h$ are witnesses for $\overline y$ and
$\overline z$ in $indec$. If $|supp \hspace{.1cm} \overline g \cap supp \hspace{.1cm} \overline h| <
\kappa$ then $\overline g$  and $\overline h$ are restrictions of $\overline f$ (meaning that apart
from a small set, their  supports are contained in $supp \hspace{.1cm} \overline f$, and on their
supports they agree with $\overline f$), and so, as $|supp \hspace{.1cm} \overline g| + |supp
\hspace{.1cm} \overline h| =  |supp \hspace{.1cm} \overline f|$, either $|supp \hspace{.1cm}
\overline g| = |supp \hspace{.1cm} \overline f|$ or $|supp \hspace{.1cm} \overline h| = |supp
\hspace{.1cm} \overline f|$, so that one  of $S_\kappa(\mu).\overline g, S_\kappa(\mu).\overline h$
is conjugate to $S_\kappa(\mu).\overline f$.  

So we suppose that $|supp \hspace{.1cm} \overline g \cap supp \hspace{.1cm} \overline h| \ge \kappa$ 
and aim for a contradiction. Since $\overline g$ and $\overline h$ commute $mod \hspace{.1cm} 
S_\kappa(\mu)$, by increasing $Y$ if necessary we may assume they commute outside $Y$. Let $Z$ be a 
typical orbit of $\langle \overline g, \overline h \rangle$ on $(supp \hspace{.1cm} \overline g 
\cap supp \hspace{.1cm} \overline h) - Y$. Then the restrictions of $\overline g$ and $\overline h$ 
to $Z$ fulfil the hypotheses of Lemma~\ref{3.5}, and so $\overline g \ast \overline h$ either has 
an orbit on $Z$ of length greater than $20$, or orbits there of length $20$ and some other length 
greater than $1$. Since this applies to all possible choices of $Z$, either there are $\ge \kappa 
\hspace{.1cm} Z$s for which there is an orbit of length greater than $20$, or there are $\ge \kappa 
\hspace{.1cm} Z$s containing an orbit of length $20$, and of some other length greater than $1$. But 
each of these is contrary to the hypothesis on $\overline f$.     $\Box$  

\vspace{.1in}

We are now able to express disjointness of certain sequences, which is the key to recovering the
appropriate ring of sets inside $S_\lambda(\mu)/S_\kappa(\mu)$. From this we shall be able to
express disjointness of involutions (meaning disjointness of their supports), which are actually
the elements we shall use to represent sets, and of more general sequences. But the first
approximation uses elements satisfying $indec$ and acting in the same way. Let us say that two such
elements $S_{\kappa}(\mu).\overline f$ and $S_{\kappa}(\mu).\overline g$ have the {\it same action}
if  $\nu_m(\overline f) \ge \kappa$ and $\nu_m(\overline g) \ge \kappa$ for the same $m > 1$.

\vspace{.1in}

\noindent $disj_1(\overline x,\overline y): indec(\overline x) \wedge indec(\overline y)
\wedge comm(\overline x,\overline y) \wedge indec(\overline x \ast \overline y)$.

\begin{lemma} For any sequences $\overline f$ and $\overline g$ of elements of $S_\lambda(\mu) -
S_\kappa(\mu)$ of length $60, S_\lambda(\mu)/S_\kappa(\mu) \models  disj_1(S_\kappa(\mu). {\overline
f},S_\kappa(\mu). {\overline g})$ if and only if $\overline f$ and $\overline g$ are indecomposable
with the same action, $|\mu \, - \, (supp \; \overline f \; \cup \; supp \; \overline g)| = \mu$, and
$|supp \; \overline f \; \cap \; supp \; \overline g| < \kappa$. \label{3.7} \end{lemma}

\noindent{\bf Proof}\hspace{.1in} It is clear that if two indecomposable sequences in 
$S_{\lambda}(\mu)/S_{\kappa}(\mu)$ have the same action, and almost disjoint supports (meaning that
the intersection of their supports has cardinality less than $\kappa$), then they commute, and their
pointwise product also is indecomposable (provided that the union of their supports does not have
small complement). Conversely suppose that the given conditions apply. Then as in the previous proof,
if the supports of $\overline f$ and $\overline g$ are not almost disjoint, then indecomposability of
$\overline f \ast \overline g$ is violated. It also follows that $\overline f, \overline g$, and
$\overline f \ast \overline g$ must all have the same action.       $\Box$

\vspace{.1in}

It is now possible to find formulae expressing the following concepts inside
$S_\lambda(\mu)/S_\kappa(\mu)$:

membership in ${\cal P}_{\lambda}(\mu)/{\cal P}_{\kappa}(\mu)$,

the boolean operations on ${\cal P}_{\lambda}(\mu)/{\cal P}_{\kappa}(\mu)$,

the action of $S_{\lambda}(\mu)/S_{\kappa}(\mu)$ on ${\cal P}_{\lambda}(\mu)/{\cal P}_{\kappa}(\mu)$.

First we represent members of ${\cal P}_{\lambda}(\mu)/{\cal P}_{\kappa}(\mu)$ by involutions, and
let $set(x)$ be the formula $x^2 = 1$, (where for present purposes it is easier to count the identity
as an `involution'). The idea is that each involution will encode its support (so for example the
identity represents the empty set). Of course this only makes any sense if we can tell when two
involutions encode the same set. 

Now let $i$ be such that $a_i$ has order 2 in $A(5)$. Then for $g \in S_\lambda(\mu)$ with $|\mu -
supp \; g| = \mu$, $S_\kappa(\mu).g$ has order $2$ if and only if there is some indecomposable
$\overline g$ with $\nu_5(\overline g) \ge \kappa$ such that $S_\kappa(\mu).g = S_\kappa(\mu).g_i$.

\noindent $disj'(x,y): set(x) \wedge set(y) \wedge \exists \overline z \exists \overline t(z_i = x
\wedge t_i =  y \wedge disj_1(\overline z,\overline t))$.

\noindent $disj(x,y): \exists x_1 \exists x_2 \exists x_3 \exists x_4 \exists y_1 \exists y_2
\exists y_3 \exists y_4(x = x_1x_2x_3x_4 \wedge y = y_1y_2y_3y_4 \wedge$

\hspace{3in} $\bigwedge_{1 \le i,j \le 4}disj'(x_i,y_j))$.

Here the idea is that $disj'$ should express disjointness of (sets encoded by) involutions, and
$disj$ should express disjointness of (the supports of) arbitrary permutations. Because of the
possibility that $\lambda = \mu^+$ we use products of four elements rather than just two, since
we need to be able to express an arbitrary group element in terms of involutions the complements of
whose supports have cardinality $\mu$.   

\noindent $subset(x,y): set(x) \wedge set(y) \wedge \forall z(disj(y,z) \rightarrow disj(x,z))$,

\noindent $sameset(x,y): set(x) \wedge set(y) \wedge \forall z(disj(y,z) \leftrightarrow disj(x,z))$,
  
\noindent $union(x,y,z): set(x) \wedge set(y) \wedge set(z) \wedge \forall t(subset(x,t) \wedge
subset(y,t) \leftrightarrow$ 

$subset(z,t))$,

\noindent $intersect(x,y,z): set(x) \wedge set(y) \wedge set(z) \wedge \forall
t(subset(t,x) \wedge subset(t,y) \leftrightarrow$ 

$subset(t,z))$,

\noindent $union_n({\overline x},y): (\forall z)(disj(z,y) \leftrightarrow \bigwedge_{i=1}^n
disj(z,x_i))$,

\noindent $map(x,y,z): set(x) \wedge set(y) \wedge sameset(z^{-1}xz,y)$,

\vspace{.1in}

The following result sums up what these formulae express.

\begin{lemma} (i) For any $f,g \in S_\lambda(\mu), S_\lambda(\mu)/S_\kappa(\mu) \models
disj'(S_\kappa(\mu).f,S_\kappa(\mu).g)$ if and only if $S_\kappa(\mu).f$ and $S_\kappa(\mu).g$ are
involutions such that $|supp \; f \, \cap \, supp \; g| < \kappa$ and $|\mu - (supp \; f \, \cup
\, supp \; g)| = \mu$. 

(ii) For any $f,g \in S_\lambda(\mu), S_\lambda(\mu)/S_\kappa(\mu) \models
disj(S_\kappa(\mu).f,S_\kappa(\mu).g)$ if and only if $|supp \; f \cap supp \; g| < \kappa$.

(iii) For any $f,g \in S_\lambda(\mu), S_\lambda(\mu)/S_\kappa(\mu) \models
subset(S_\kappa(\mu).f,S_\kappa(\mu).g)$ if and only if $S_\kappa(\mu).f$ and $S_\kappa(\mu).g$ are
involutions such that $|supp \; f - supp \; g| < \kappa$. 

(iv) For any $f,g \in S_\lambda(\mu), S_\lambda(\mu)/S_\kappa(\mu) \models
sameset(S_\kappa(\mu).f,S_\kappa(\mu).g)$ if and only if $S_\kappa(\mu).f$ and $S_\kappa(\mu).g$ are
involutions such that $|supp \; f - supp \; g|$, $|supp \; g - supp \; f| < \kappa$.

(v) For any $f,g,h \in S_\lambda(\mu), S_\lambda(\mu)/S_\kappa(\mu) \models
union(S_\kappa(\mu).f,S_\kappa(\mu).g,S_\kappa(\mu).h)$ if and only if $S_\kappa(\mu).f,
S_\kappa(\mu).g$, and $S_\kappa(\mu).g$ are involutions such that
$supp \; f  \cup supp \; g$ and $supp \; h$ differ by a set of cardinality $< \kappa$.

(vi) Similarly for intersections.

(vii) For any ${\overline f}, g \in S_\lambda(\mu), S_\lambda(\mu)/S_\kappa(\mu) \models
union_n(S_\kappa(\mu).{\overline f},S_\kappa(\mu).g)$ if and only if $\bigcup_{i=1}^n supp(f_i)$ 
and $supp \;g $ differ by a set of cardinality $< \kappa$.

(viii) For any $f,g,h \in S_\lambda(\mu), S_\lambda(\mu)/S_\kappa(\mu) \models
map(S_\kappa(\mu).f,S_\kappa(\mu).g,S_\kappa(\mu).h)$ if and only if $S_\kappa(\mu).f$ and
$S_\kappa(\mu).g$ are involutions and $(supp \; f)h$ and $supp \; g$ differ by a set of cardinality $<
\kappa$.
                           \label{3.8}  \end{lemma}

\noindent{\bf Proof}\hspace{.1in} (ii) follows from the fact that any permutation may be written as a 
product of two involutions, and any involution may be written as a product of two involutions
the complement of whose support has cardinality $\mu$. The rest of the proof is straightforward.    
$\Box$

\begin{corollary} The ring of sets ${\cal P}_\lambda(\mu)/{\cal P}_\kappa(\mu)$ and the natural
action of $S_\lambda(\mu)/S_\kappa(\mu)$ on this ring are interpretable inside the group
$S_\lambda(\mu)/S_\kappa(\mu)$.  \label{3.9}  \end{corollary}

This result is due to Rubin \cite{Rubin1} Theorem 4.3, but using different methods. An alternative
route to the same conclusion, avoiding so much detail on permutation representations, starts by
interpreting ${\cal P}_\lambda(\mu)/{\cal P}_\kappa(\mu)$ in $S_\lambda(\mu)/S_\kappa(\mu)$ using
parameters ${\overline {f^*}},{\overline {f^{**}}}$. The first of these acts as $A(5)$ with orbits of
degree $5$ and $1$ only, and with the aid of the second, disjointness can be expressed more rapidly.
The parameters are then eliminated at a later stage.

\section{Interpreting ${\cal M}_{\kappa \lambda \mu}$ in $S_\lambda(\mu)/S_\kappa(\mu)$}

This is carried out as follows:
 
members of $IS_n$ are represented by `pure' $n$-tuples, being those for which almost all orbits are
isomorphic, modulo isomorphism of this action,

members of $Card^-$ are represented by group elements which encode sets, (that is, involutions),
modulo the relation of having equal cardinality,

members of $F_n$ are represented by $n$-tuples of group elements, modulo conjugacy.

In addition we have to show definability of the relations and functions in the signature.

First we show how to distinguish the case $\lambda = \mu^+$ (which has already required special
treatment in the previous section). We use the formula

$$max: (\exists x)(\forall y)(disj(x,y) \rightarrow y = 1)$$

\noindent (expressing that $\lambda$ has its {\em maximum} value)

\begin{remark} $\lambda = \mu^+$ if and only $S_\lambda(\mu)/S_\kappa(\mu) \models max$.     
                                       \label{4.1}      \end{remark}

To carry out the interpretation more formally we require the following formulae:

\noindent $disj_n(\overline x,\overline y): \bigwedge_{1 \le i,j \le n}disj(x_i,y_j)$,

\noindent $restr_n(\overline x,\overline y): \exists \overline z(disj_n(\overline x,\overline z)
\wedge \overline x \ast \overline z = \overline y)$,

\noindent $\overline x = 1: \bigwedge_{i=1}^n x_i = 1$,

\noindent $compat_n(\overline x,\overline y): \exists \overline z \exists t(\overline z \neq 1
\wedge restr_n(\overline z,\overline x) \wedge restr_n(\overline z^t,\overline y))$,

\noindent $pure_n(\overline x): \forall \overline y \forall \overline z (\overline y \neq 1 \wedge
\overline z \neq 1 \wedge restr_n(\overline y,\overline x) \wedge restr_n(\overline z,\overline x)$

$\rightarrow compat_n(\overline y,\overline z)) \wedge (\neg max \rightarrow {\overline x} \neq 1)$,

\noindent $iso_n(\overline x,\overline y): pure_n(\overline x) \wedge pure_n(\overline y) \wedge
(compat_n(\overline x,\overline y) \vee {\overline x} = {\overline y} = 1)$.

\begin{lemma} (i) For any finite sequence $\overline f$ of members of $S_\lambda(\mu),
S_\lambda(\mu)/S_\kappa(\mu) \models pure_n(S_\kappa(\mu).\overline f)$ if and only if the
non-trivial actions of $\overline f$ on all but a small union of the orbits of $\langle \overline f
\rangle$ are isomorphic or, if $\lambda = \mu^+$, almost all orbits have size $1$..

(ii) For any sequences $\overline f, \overline g$ in $S_\lambda(\mu), S_\lambda(\mu)/S_\kappa(\mu)
\models iso_n(S_\kappa(\mu).{\overline f},S_\kappa(\mu).{\overline g})$ if and only if the actions of
$\overline f$ and $\overline g$ on all but a small union of orbits of $\langle \overline f \rangle,
\langle \overline g \rangle$ have the same isomorphism type in $IS_n$.  \label{4.2}  \end{lemma}

Note that it is not enough to talk of the actions of $\langle \overline f \rangle$ on its orbits; we
need to distinguish the generating tuple $\overline f$ in order to capture $IS_n$. Observe that the
final parts of the formulae $pure_n$ and $iso_n$ cover the case $\lambda = \mu^+$, and correspond to
the remark in parentheses in Definition \ref{2.1}(ii). Similar remarks apply to the treatment of
$F_n$.

As mentioned above, for sort $2$ we just use involutions, this time modulo the equivalence relation
given by

\noindent $samecard(x,y): set(x) \wedge set(y) \wedge \exists x_1 \exists x_2 \exists y_1 \exists
y_2 (disj(x_1,x_2) \wedge disj(y_1,y_2) \wedge$ 

$x = x_1x_2 \wedge y = y_1 y_2 \wedge conj(x_1,y_1) \wedge conj(x_2,y_2))$.

This is slightly more complicated than the expected `$set(x) \wedge set(y) \wedge conj(x,y)$' in
view of the case $\lambda = \mu^+$. And the sorts $3$ have already been remarked on.

It remains to show that the relations and functions of ${\cal M}_{\kappa \lambda \mu}$ are definable.

First the ordering $\le$ (and hence $<$) on $Card$ is definable by

\noindent $lesseq(x,y): set(x) \wedge set(y) \wedge (\exists z)(subset(z,y) \wedge samecard(x,z))$.

To define $Eq^1,Eq,Prod^1$ and $Prod$ we use 

\noindent $eq^1(x_1,x_2): pure_2(x_1,x_2) \wedge$ $x_1 = x_2, eq(x_1,x_2): x_1 = x_2$, 

\noindent $prod^1(x_1,x_2,x_3): pure_3(x_1,x_2,x_3) \wedge x_1x_2 = x_3$ 

and

\noindent $prod(x_1,x_2,x_3): x_1x_2 = x_3$ 

\noindent respectively. Note that there is a slight difference between $Eq^1$ and $Eq$ (and between
$Prod^1$ and $Prod$), since in the former case $(x_1,x_2)$ is meant to represent a member of $IS_2$,
but in the latter, of $F_2$.

We may define $Proj_n^1$ and $Proj_n$ by

\noindent $proj_n^1((x_1,\ldots,x_{n+1}),(y_1,\ldots,y_n)): pure_{n+1}(x_1,\ldots,x_{n+1}) \wedge 
pure_n(y_1,\ldots,y_n) \wedge$

$iso_n(x_1, \ldots , x_n, y_1,\ldots,y_n)$, and

\noindent $proj_n((x_1,\ldots,x_{n+1}),(y_1,\ldots,y_n)):
conj_n((x_1,\ldots,x_n),(y_1,\ldots,y_n))$,

\noindent and $App_n$ by

\noindent $app_n({\overline x},{\overline y},z): pure_n({\overline y}) \wedge [((\exists {\overline
t}) (pure_n({\overline t}) \wedge compat_n({\overline y},{\overline t}) \wedge 
restr_n({\overline t},{\overline x})$

$\wedge (\forall {\overline u})(restr_n({\overline t},{\overline u}) \wedge restr_n({\overline
u},{\overline x}) \wedge pure_n({\overline u} ) \rightarrow {\overline t} = {\overline u}) \wedge 
(\exists v)(union_n({\overline t} ,v)$

$\wedge \; samecard(v,z))) \vee ((\forall {\overline t})(compat_n({\overline y}, {\overline t})
\rightarrow \neg restr_n({\overline t},{\overline x})) \wedge z = 1)]$,

\noindent which we may paraphrase as `either there is a maximal pure restriction $\overline t$ of
$\overline x$ compatible with $\overline y$ and of cardinality (coded by) $z$, or $\overline x$ has 
no restriction compatible with $\overline y$ and $z = 1$ (that is, codes $0$)'. If $\lambda =
\mu^+$, $app_n$ is modified to cover the case ${\overline y} = 1$, and if $\kappa = \aleph_0$ we have
to count orbits rather than their union, and the statement about $v$ is modified to express `there is
a set having the same cardinality as $z$ which intersects each orbit of $\overline t$ and is minimal
subject to this'. To justify this we further note that the case $\kappa = \aleph_0$ can be
distinguished by the sentence

$$(\exists x)(\forall y)(restr_1(y,x) \rightarrow (y = 1 \vee y = x)).$$

We have proved the following:

\begin{theorem} ${\cal M}_{\kappa \lambda \mu}$ is interpretable in the group
$S_\lambda(\mu)/S_\kappa(\mu)$. \label{4.3} \end{theorem}

We remark that `interpretability' here is taken in the usual sense (called `explicit interpretability'
in \cite{Shelah1}). This means for instance that, rather than just transferring the first order
properties, we are able to deduce that whenever
$S_{\lambda_1}({\mu}_1)/S_{\kappa_1}(\mu_1) \cong S_{\lambda_2}({\mu}_2)/S_{\kappa_2}(\mu_2)$ then
${\cal M}_{\kappa_1 \lambda_1 \mu_1} \cong {\cal M}_{\kappa_2 \lambda_2 \mu_2}$, and hence to try to
distinguish the groups $S_\lambda(\mu)/S_\kappa(\mu)$ up to isomorphism as well as up to elementary
equivalence. But for us here the following is the point.

\begin{corollary} ${\cal M}_{\kappa_1 \lambda_1 \mu_1} \equiv {\cal M}_{\kappa_2 \lambda_2 \mu_2}$ if
and only if $S_{\lambda_1}(\mu_1)/S_{\kappa_1}(\mu_1) \equiv$

\noindent $S_{\lambda_2}(\mu_2)/S_{\kappa_2}(\mu_2)$.  \label{4.4}    \end{corollary}

\noindent{\bf Proof}\hspace{.1in} This follows from Theorems \ref{2.6} and \ref{4.3}.
$\Box$

\vspace{.1in}

In the next sections we give more details about the circumstances under which ${\cal M}_{\kappa_1
\lambda_1 \mu_1} \equiv {\cal M}_{\kappa_2\lambda_2\mu_2}$.

\section{Refinements and the case $cf(\kappa) > 2^{\aleph_0}$}

We now make some remarks about distinguishing the elementary theories of
$S_\lambda(\mu)/S_\kappa(\mu)$ for different values of $\kappa, \lambda, \mu$, which by Corollary
\ref{4.4} is equivalent to distinguishing the ${\cal M}_{\kappa \lambda \mu}$. In the first
place, according to Remark \ref{4.1}, the case $\lambda = \mu^+$ can be singled out in
$S_\lambda(\mu)/S_\kappa(\mu)$ by means of the sentence $max$ of the language of group theory, and
hence also by a suitable sentence in ${\cal M}_{\kappa \lambda \mu}$. So we may treat the cases
$\lambda \le \mu, \lambda = \mu^+$ separately. Now when $\lambda \le \mu$ the cardinal $\mu$ actually
plays no part at all in the structure ${\cal M}_{\kappa \lambda \mu}$, so we at once see that for
fixed $\kappa \le \lambda$, all the ${\cal M}_{\kappa \lambda \mu}$ with $\mu \ge \lambda$ are
elementarily equivalent. More is even true at this stage, since many of the ${\cal M}_{\kappa \lambda
\mu}$ are in fact isomorphic. For instance if $cf(\kappa_1), cf(\kappa_2) \ge (2^{\aleph_0})^+$ then
${\cal M}_{\kappa_1 \kappa_1^+ \mu_1} \cong {\cal M}_{\kappa_2 \kappa_2^+ \mu_2}$ (for $\mu_1 \ge
\kappa_1^+, \mu_2 \ge \kappa_2^+$) since in this case $Card^- = \{ 0, \kappa_1 \}, \{ 0, \kappa_2 \}$
respectively, and similarly ${\cal M}_{\kappa_1 \kappa_1^{++} \mu_1} \cong {\cal M}_{\kappa_2
\kappa_2^{++} \mu_2}$ etc.

We know of course that $Th({\cal M}_{\kappa \lambda \mu})$ can only take at most $2^{\aleph_0}$
values, and so there will be many pairs of distinct triples giving elementarily equivalent models.
In \cite{Shelah2} this was however illustrated more explicitly, and we carry out a similar analysis
here. There a characterization of elementary equivalence was provided based on the second order
theory of certain many-sorted ordinal structures, whose sorts all had cardinality $\le 2^{\aleph_0}$,
and we give a parallel treatment. While doing so we give a few more details about the material from
\cite{Shelah2} (which in its turn is related to \cite{Kino}). First we  show how a suitable second
order logic can be represented in the structures ${\cal M}_{\kappa \lambda \mu}$. Small modifications
are made in the case $\kappa = \aleph_0$ (distinguishable in the language of group theory), which we
do not spell out explicitly.

To represent subsets of $IS_n$ in ${\cal M}_{\kappa \lambda \mu}$ is rather straightforward, but 
subsets of $IS_m \times IS_n$ are harder to deal with. We use `products' (similar to the method of
section 3) to help us to do this. We say that $t \in IS_{m+n}$ is a {\em product} of
$t_1 \in IS_m$ and $t_2 \in IS_n$ if $t$ has the form $((A,{\overline f}))_{\cong}$ where $A =
\{\alpha_{xy}: x \in X, y \in Y\}$, the $\langle f_1, \ldots, f_m \rangle$-orbits of $A$ are
$\{\alpha_{xy}: y \in Y\}$ for $x \in X$, all of type $t_1$, the $\langle f_{m+1}, \ldots, f_{m+n}
\rangle$-orbits of $A$ are $\{\alpha_{xy}: x \in X\}$ for $y \in Y$, all of type $t_2$, and the actions
of $f_i$ and $f_j$ on $A$ for $1 \le i \le m < j \le m+n$ commute. We say that $h \in F_{m+n}$ is a
{\em product} if whenever $h = (h')_{{\cal E}_{m+n}}$, $\sum \{h'(t): t \in IS_{m+n} \wedge t$ not a 
product$\} < \kappa$. 

The idea here is that if $\overline f$ acts as a product on almost all of its orbits, then we can
uniquely recover its actions on the first $m$ and last $n$ co-ordinates, so that products provide a
way of encoding sets of ordered pairs. As illustrated in section 3 however, the actions of tuples may
commute without their being a product, and so the natural condition to try to capture expressibility
as a product, namely commutativity, does not work. This time however this does not matter; the point
being that when two actions commute, and together generate a transitive action, the projections onto
the two sets of co-ordinates are uniquely determined. Let us therefore say that $h \in F_{m+n}$ is a
{\em product} if whenever $h = (h')_{{\cal E}_{m+n}}$, 
$$\sum \{h'(t): t \in IS_{m+n} \wedge t = ((A,({\overline f}_1,{\overline f}_2)))_{\cong}
\rightarrow \neg comm_{m,n}({\overline f}_1,{\overline f}_2)\} < \kappa.$$ 

We now represent subsets of $IS_n$ of cardinality $< \lambda$ by $h \in F_n$ such that $(\forall t \in
IS_n)h(t) \le \kappa$, and subsets of $IS_m \times IS_n$ of cardinality $< \lambda$ by $k \in F_{m+n}$
which are products and such that $(\forall t \in IS_{m+n})k(t) \le \kappa$. The subset of $IS_n$
encoded by $h$ is then $\{t: h(t) = \kappa\}$ and the subset of $IS_m \times IS_n$ encoded by $k$ is
$\{(t_1,t_2) \in IS_m \times IS_n: (\exists t \in IS_{m+n})(k(t) = \kappa$ and $t_1,t_2$ are the
members of $IS_m,IS_n$ determined on co-ordinates 1 to $m$ and $m+1$ to $m+n$ respectively)$\}$.
(The definition of `product' ensures that these are uniquely determined from $t$, since if the
actions of
${\overline f}_1$ and ${\overline f}_2$ on the $\overline f$-orbit $A$ commute then they each
preserve the set of orbits of the other, and all actions of ${\overline f}_1$ on its orbits are
isomorphic, and similary for ${\overline f}_2$.)

\begin{lemma}  There are formulae of the language of ${\cal M}_{\kappa \lambda \mu}$ expressing
the following:

(i) $h$ encodes a subset of $IS_n$,

(ii) $k$ encodes a subset of $IS_m \times IS_n$,

(iii) $mem_n(t,h): t$ lies in the set encoded by $h$,

$mem_{m,n}(t_1,t_2,k): (t_1,t_2)$ lies in the set encoded by $k$,

(iv) $equal_n(h,h'): h,h'$ encode the same subset of $IS_n$,

$equal_{m,n}(k,k'): k,k'$ encode the same subset of $IS_m \times IS_n$,

(v) $fun_{m,n}(k): k$ encodes a function (from a subset of $IS_m$ into $IS_n$),

(vi) $one$-$onefun_{m,n}(k): k$ encodes a 1--1 function.   \label{5.1}     \end{lemma}

\noindent {\bf Proof} \hspace{.1in} (i) $h$ encodes a subset of $IS_n$ if and only if $(\forall t
\in IS_n)App_n(h,t) \le \kappa$.

(ii) By appeal to Theorem \ref{2.6} we may express projections of $h \in F_{m+n}$ to co-ordinates 1
to $m$ and $m+1$ to $m+n$, and then use the formula $comm_{m,n}$. 

(iii) $mem_n(t,h)$ is taken as $App_n(h,t) = \kappa$.

For $mem_{m,n}(t_1,t_2,k)$ we take $(\exists t \in IS_{m+n})(App_{m+n}(h,t) = \kappa \wedge t_1, t_2$
are the projections of $t$ onto co-ordinates $1$ to $m$ and $m+1$ to $m+n$ respectively). (The fact
that we can express these more generalized projections here follows by appeal to Theorem \ref{2.6},
though they could also have been included in the signature of the ${\cal M}_{\kappa \lambda \mu}$ if
desired.)    

(iv), (v), and (vi) follow from (iii).      $\Box$

\vspace{.1in}   

For the remainder of this section we specialize to the case $cf(\kappa) > 2^{\aleph_0}$, to avoid
complications. We return to the general case in section 6. One of the benefits of assuming
$cf(\kappa) > 2^{\aleph_0}$ is that we can dispense altogether with the equivalence relations
${\cal E}_n$. This is because any function from $F_n$ to $Card$ is ${\cal E}_n$-equivalent to a unique
function from $F_n$ into $Card^-$ (obtained by replacing all values below $\kappa$ by $0$). Various
other simplifications and interdefinabilities in this case are described in the following theorem.

\begin{theorem} Suppose that $\kappa < \lambda \le \mu^+$ and $cf(\kappa) > 2^{\aleph_0}$. Then 

(i) for each $h \in F_n$ there is a unique ${\cal E}_n$-representative which is a function from
$IS_n$ to $Card^-$ (so that from now on in this section we dispense with ${\cal E}_n$ and regard $F_n$
as a subset of $(Card^-)^{IS_n}$),

(ii) for each $n, Sum_n: F_n \rightarrow Card^-$ given by $Sum_n(h) = \sum \{ h(t): t \in IS_n
\}$ is definable in ${\cal M}_{\kappa \lambda \mu}$,

(iii) $Eq, Prod$, and $Proj_n$ are all definable in 
$${\cal M}^*_{\kappa \lambda \mu} = ((IS_n)_{n \ge 1}, Card^-, (F_n)_{n \ge 1}; Eq^1, Prod^1,
(Proj_n^1)_{n \ge 1}, <, (App_n)_{n \ge 1}),$$ 
and conversely $Eq^1, Prod^1$, and $Proj_n^1$ are all definable in  
$$((IS_n)_{n \ge 1}, Card^-, (F_n)_{n \ge 1}; <, (App_n)_{n \ge 1}, Eq, Prod, (Proj_n)_{n \ge 1}).$$ 
                  \label{5.2}   \end{theorem}

\vspace{-.1in}
\noindent{\bf Proof}\hspace{.1in}(ii) This is because $Sum_n(h)$ may also be written as $\sup 
\{ h(t): t \in IS_n \}$ (since $\kappa > 2^{\aleph_0}$), so that
$$Sum_n(h) = \alpha \Leftrightarrow (\forall t)App_n(h,t) \le \alpha \: \wedge \: (\forall \beta <
\alpha)(\exists t)(\beta < App_n(h,t)).$$

\noindent (iii) $Eq = \{ h \in F_2: (\forall t \in IS_2)(App_2(h,t) \neq 0 \rightarrow Eq^1(t)) \}$,

$Prod = \{ h \in F_3: (\forall t \in IS_3)(App_3(h,t) \neq 0 \rightarrow Prod^1(t)) \}$.

\noindent For $Proj_n$ we remark that $B_{t't} \neq \emptyset \Leftrightarrow Proj_n^1(t',t)$, and so
\vspace{-.1in} 
$$Proj_n(h)(t) = \sum \{ h(t'): Proj_n^1(t',t) \} = \sup \{ h(t'): Proj_n^1(t',t) \}.$$
\vspace{-.2in}

\noindent As in (ii) we see that 
$$\begin{array}{c}
Proj_n(h)(t) = \alpha \Leftrightarrow (\forall t')(Proj_n^1(t',t) \rightarrow App_{n+1}(h,t') \le
\alpha) \\
\wedge \: (\forall \beta < \alpha)(\exists t' \in IS_{n+1})(Proj_n^1(t',t) \: \wedge \: \beta <
App_{n+1}(h,t')).
     \end{array} $$
Conversely we have
$$\begin{array}{c}
Eq^1 = \{ t \in IS_2: (\exists h \in F_2)(Eq(h) \: \wedge \: App_2(h,t) \neq 0) \}, \\
Prod^1 = \{ t \in IS_3: (\exists h \in F_3)(Prod(h) \: \wedge \: App_3(h,t) \neq 0) \},
     \end{array} $$
and
$$\;\;\;Proj_n^1(t_1,t_2) \Leftrightarrow (\forall h \in F_{n+1})(App_{n+1}(h,t_1) \neq 0 \rightarrow
Proj_n(h)(t_2) \neq 0). \;\;\; \Box$$   
 
\vspace{.1in}

This theorem tells us that when $cf(\kappa) > 2^{\aleph_0}$ it suffices to consider the structures
${\cal M}^*_{\kappa \lambda \mu}$, and here, since the members of $F_n$ are now viewed as functions
from $IS_n$ to $Card^-$, this amounts to a version of second order logic on the sorts $1$ and $2$,
together with $Eq^1, Prod^1, Proj_n^1$, and $<$. The sorts $IS_n$ are independent of $\kappa,
\lambda, \mu$, and so the main point is to analyse $Card^-$. We give an analysis of this situation
similar to that in \cite{Shelah2} which involves defining suitable `small' ordinals (meaning
of cardinality $\le 2^{\aleph_0}$), sufficient to capture the elementary theory. 

In what follows we extend the definition of `cofinality' to zero or successor ordinals by
letting $cf(0) = 0$ and $cf(\alpha + 1) = 1$. Let $\Omega = (2^{\aleph_0})^+$. Then any ordinal
$\alpha$ may be written uniquely in the form  $$\alpha = \Omega^\omega.\alpha_\omega + \ldots +
\Omega^n.\alpha_{[n]} + \ldots + \Omega.\alpha_{[1]} + \alpha_{[0]}$$ where $\Omega^\omega$ is the
ordinal power, $\alpha_{[n]} < \Omega$ for $n \in \omega$, and $\{ n: \alpha_{[n]} \neq 0 \}$ is
finite. We write $\alpha [n] = \Omega^\omega.\alpha_\omega + \ldots + \Omega^{n+1}.\alpha_{[n+1]}$
and let  
\[ \alpha^{[n]} = \left\{ \begin{array}{l} 
        1 + cf(\alpha [n]) \hspace{.2in} \mbox{ if $cf(\alpha [n]) < \Omega$,}            \\
        0   \hspace{1.2in}   \mbox{ otherwise}
\end{array} \right . \]
For ordinals $\alpha, \beta$, and $k \in \omega$, let $\alpha \sim_k \beta$ if $\alpha_{[l]} =
\beta_{[l]}$ and $\alpha^{[l]} = \beta^{[l]}$ for all $l \le k$, and $\alpha \sim \beta$ if $\alpha
\sim_k \beta$ for all $k$. For a set of ordinals $\{ \alpha \} \cup A$ let $\gamma (\alpha,A)$ be the
order-type of $\{ \beta \le \alpha: (\forall \gamma \in A)(\gamma < \alpha \rightarrow \gamma <
\beta) \}$. Note that this set is the final segment of $\alpha \cup \{ \alpha \}$ consisting of all
(strict) upper bounds of $A \cap \alpha$. 

\vspace{.1in}

The following lemma is stated in \cite{Shelah1} and a proof outlined. The result is also related to
\cite{Kino}. We give fuller details here for the reader's benefit.

\begin{lemma} (i) $\sim_k$ is an equivalence relation. For each $\alpha$ there is $\beta <
\Omega^{k+2}$ with $\alpha \sim_k \beta$; and if $\alpha \sim_k \beta$ then $\alpha < \Omega^{k+1}$
if and only if $\beta < \Omega^{k+1}$, and each of these implies $\alpha = \beta$.

(ii) If $\alpha \ge \Omega^{k+1}$ then for any $\beta, \alpha \sim_k \beta + \alpha$.

(iii) If $\alpha_\gamma \sim_k \beta_\gamma$ for each $\gamma < \delta$, then $\sum_{\gamma <
\delta}\alpha_\gamma \sim_k \sum_{\gamma < \delta}\beta_\gamma$.

(iv) If $\alpha \sim_{k+1} \beta$, and $A \subseteq \alpha$ with $|A| < \Omega$, there is an
order-preserving map $F: A \rightarrow \beta$ such that for each $a \in A \cup \{ \alpha \},
\gamma(a,A) \sim_k \gamma (F(a),F(A))$ (where $F(\alpha)$ is taken to equal $\beta$).
   \label{5.3}     \end{lemma}

\noindent{\bf Proof}\hspace{.1in} (i) Let $\beta = \Omega^{k+1}.\beta_{[k+1]} + \Omega^k.\alpha_{[k]} + \ldots
+ \Omega.\alpha_{[1]} + \alpha_{[0]}$ where $\beta_{[k+1]}$ is given as follows:
\[
\beta_{[k+1]} = \left\{ \begin{array}{l} 
      \alpha_{[n]} \hspace{.6in} \mbox{ if  $\alpha_{[n]} \neq 0$ for some least $n > k$,}    \\
    \omega \hspace{.72in} \mbox{ if $\alpha_{[n]} = 0$ for all $n > k$ and $\alpha_\omega$ is a
                                         successor,} \\
    cf(\alpha_\omega) \hspace{.4in} \mbox{ if  $\alpha_{[n]} = 0$ for all $n > k, \alpha_\omega$
                                       a limit ordinal and }  \\
          \hspace{3in}      cf(\alpha_\omega) < \Omega \\
        0        \hspace{1.2in}     \mbox{    otherwise}
\end{array} \right . \]
Then $\alpha_{[l]} = \beta_{[l]}$ for $l \le k$ is immediate. If $\alpha < \Omega^{k+1}$ the final
clause applies, so $\beta = \alpha < \Omega^{k+1}$. Also if $\beta < \Omega^{k+1}$ then $\alpha <
\Omega^{k+1}$ so the last part also follows.

Now suppose the first clause applies. Then if $l \le k$,
$$\alpha [l] = \Omega^\omega.\alpha_\omega + \ldots + \Omega^n.\alpha_{[n]} + \Omega^k.\alpha_{[k]}
+ \ldots + \Omega^{l+1}.\alpha_{[l+1]}$$
$$\mbox{and } \beta [l] = \Omega^{k+1}.\alpha_{[n]} + \Omega^k.\alpha_{[k]}
+ \ldots + \Omega^{l+1}.\alpha_{[l+1]}$$
which have equal cofinalities as $\alpha_{[n]} \neq 0$. If the second or third clause applies, then
$$\alpha [l] = \Omega^\omega.\alpha_\omega + \Omega^k.\alpha_{[k]} + \ldots +
                                                      \Omega^{l+1}.\alpha_{[l+1]}$$ 
$$\mbox{and } \beta [l] = \Omega^{k+1}.\beta_{[k+1]} + \Omega^k.\alpha_{[k]} + \ldots +
                                                      \Omega^{l+1}.\alpha_{[l+1]}.$$
If $\alpha^{[l]} \neq \beta^{[l]}$ then $\alpha_{[k]} = \ldots = \alpha_{[l+1]} = 0$ so $\alpha [l]
= \Omega^\omega.\alpha_\omega$ and $\beta [l] = \Omega^{k+1}.\beta_{[k+1]}$. But if clause $2$
applies, $cf(\alpha [l]) = \omega = cf(\beta [l])$, and if clause $3$ applies, $cf(\alpha [l]) =
cf(\alpha_\omega) = \beta_{[k+1]} = cf(\beta [l])$ after all.

(ii) As $\alpha \ge \Omega^{k+1}, \alpha_\omega \neq 0$, or $\alpha_{[n]} \neq 0$ for some $n > k$.
Write $\alpha = \Omega^N.\alpha_{[N]} + \ldots + \alpha_{[0]}$ where $N > k, \alpha_{[N]} \neq 0$
(and where $N = \omega, \alpha_{[N]} = \alpha_\omega$ is allowed). Writing $\beta$ in a similar
way, if $n < N, \Omega^n.\beta_{[n]} + \Omega^N.\alpha_{[N]} = \Omega^N.\alpha_{[N]}$, and so
$\beta + \alpha = \Omega^\omega.\beta_\omega + \ldots + \Omega^N.(\beta_{[N]} + \alpha_{[N]}) +
\Omega^{N-1}.\alpha_{[N-1]} + \ldots + \alpha_{[0]}$. For $l \le k$ we have
$$\alpha [l] = \Omega^N.\alpha_{[N]} + \ldots + \Omega^{l+1}.\alpha_{[l+1]}$$
$$\mbox{and } (\beta + \alpha)[l] = \Omega^\omega.\beta_\omega + \ldots + \Omega^N.(\beta_{[N]} +
\alpha_{[N]}) + \ldots + \Omega^{l+1}.\alpha_{[l+1]}.$$
The only way in which $cf(\alpha [l])$ can be unequal to $cf((\beta + \alpha)[l])$ is for $\alpha
[l] = \Omega^N.\alpha_{[N]} \neq \Omega^\omega.\beta_\omega + \ldots + \Omega^N.(\beta_{[N]} +
\alpha_{[N]}) = (\beta + \alpha)[l]$. But if $\alpha_{[N]}$ is a limit ordinal,
$cf(\Omega^N.\alpha_{[N]}) = cf(\alpha_{[N]}) = cf(\Omega^N.(\beta_{[N]} + \alpha_{[N]}))$, and if
it is a successor, both cofinalities are equal to $cf(\Omega^N)$, so we deduce that $\alpha^{[l]} =
\beta^{[l]}$, and hence that $\alpha \sim_k \beta$.

(iii) If $\delta = 0$ or $1$ the result is immediate. Next suppose $\delta = 2$. If $\alpha_1 \ge
\Omega^{k+1}$ then by (i) also $\beta_1 \ge \Omega^{k+1}$ (and vice versa), so by (ii) $\alpha_0 +
\alpha_1 \sim_k \alpha_1 \sim_k \beta_1 \sim_k \beta_0 + \beta_1$. Otherwise if $m$ is greatest such
that $\alpha_{1[m]} \neq 0$ then $m \le k$ and also $m$ is the greatest such that $\beta_{1[m]} \neq
0$, and $\alpha_0 + \alpha_1 =$  $$\Omega^\omega.\alpha_{0 \omega} + \ldots +
\Omega^{k+1}.\alpha_{0[k+1]} + \ldots + \Omega^m.(\alpha_{0[m]} + \alpha_{1 [m]}) +
\Omega^{m-1}.\alpha_{1 [m-1]} + \ldots + \alpha_{1 [0]},$$  
and $\beta_0 + \beta_1 =$  
$$\Omega^\omega.\beta_{0 \omega} + \ldots + \Omega^{k+1}.\beta_{0 [k+1]} + \ldots + \Omega^m.(\beta_{0
[m]} + \beta_{1 [m]}) + \Omega^{m-1}.\beta_{1 [m-1]} + \ldots + \beta_{1 [0]}.$$ 
>From $\alpha_i^{[l]} = \beta_i^{[l]}$, for $i = 0,1, l \le k$ it follows that $(\alpha_0 +
\alpha_1)^{[l]} = (\beta_0 + \beta_1)^{[l]}$.

We now prove the general case by transfinite induction. The successor case follows easily from the
case $\delta = 2$. Suppose therefore that $\delta$ is a limit ordinal. Since $\alpha_\gamma \sim_k
\beta_\gamma$ for $\gamma < \delta, \alpha_\gamma = 0 \Leftrightarrow \beta_\gamma = 0$, so we
ignore any zero terms. Thus $cf(\sum_{\gamma < \delta}\alpha_\gamma) = cf( \sum_{\gamma <
\delta}\beta_\gamma)$ ($= cf(\delta$)). Also by (i), $\alpha_\gamma \ge \Omega^{k+1} \Leftrightarrow
\beta_\gamma \ge \Omega^{k+1}$, and so $\{ \gamma < \delta: \alpha_\gamma \ge \Omega^{k+1} \}$ is
unbounded $\Leftrightarrow \{ \gamma < \delta: \beta_\gamma \ge \Omega^{k+1} \}$ is unbounded. If
each of these is unbounded, $\sum_{\gamma < \delta}\alpha_\gamma = \Omega^{k+1}.\alpha^*, 
\sum_{\gamma < \delta}\beta_\gamma = \Omega^{k+1}.\beta^*$ for some $\alpha^*, \beta^*$. Otherwise
for some $\gamma_0 < \delta, (\forall \gamma \ge \gamma_0)(\alpha_\gamma, \beta_\gamma <
\Omega^{k+1})$ and as $\alpha_\gamma \sim_k \beta_\gamma$, by (i) $(\forall \gamma \ge
\gamma_0)(\alpha_\gamma = \beta_\gamma)$. In each case it follows that $\sum_{\gamma <
\delta}\alpha_\gamma \sim_k \sum_{\gamma < \delta}\beta_\gamma$.

(iv) Given $\alpha \sim_{k+1} \beta$ and $A \subseteq \alpha, |A| < \Omega$ we write $\alpha =
\alpha' + \xi, \beta = \beta' + \xi$ where $\alpha', \beta'$ are divisible by $\Omega^{k+2}$, and
$cf(\alpha') = cf(\beta')$ or $cf(\alpha'), cf(\beta') \ge \Omega$. First suppose $\xi = 0$.

If $cf(\alpha), cf(\beta) \ge \Omega$ we define $F: A \cup \{ \alpha \} \rightarrow \beta \cup \{
\beta \}$ by $F(\alpha) = \beta$ and otherwise inductively so that for each $a \in A, \gamma(a,A)
\sim_k \gamma(F(a),F(A))$ and $\gamma(F(a), F(A)) < \Omega^{k+2}$. Suppose that $F(a')$ has been
defined for $a' < a$ having these properties. Then $\gamma(a,A)$ is known and $F(a)$ has to be
chosen. This is possible by (i), and as $|A| < \Omega$ and $\Omega$ is regular, $F(a) <
\Omega^{k+2} \le \beta$. Moreover $\gamma(\alpha,A) \sim_k \gamma(\beta,F(A))$ is clear (since each
of these order-types is cofinal with a positive multiple of $\Omega^{k+2}$).

Next if $cf(\alpha) = cf(\beta) < \Omega$ we may write $\alpha = \sum_{\gamma <
\lambda}\alpha_\gamma, \beta =  \sum_{\gamma < \lambda}\beta_\gamma$ where each $\alpha_\gamma,
\beta_\gamma$ has cofinality $\ge \Omega$ and is divisible by $\Omega^{k+2}$. By the first case we
define $F: A \cap \{ \xi: \sum_{\delta < \gamma}\alpha_\delta \le \xi < \sum_{\delta \le
\gamma} \alpha_\delta \} \rightarrow \{ \xi: \sum_{\delta < \gamma}\beta_\delta \le \xi < \sum_{\delta
\le \gamma} \beta_\delta \}$ for each $\gamma < \lambda$ and put the pieces together.

Finally for the case $\xi \neq 0$ we define $F: A \cap \alpha' \rightarrow \beta'$ as above and let
$F(\alpha' + \gamma) = \beta' + \gamma$ for $\gamma < \xi$ and $\alpha' + \gamma \in A$.
   $\Box$

\vspace{.1in}

Now we can prove the required bi-interpretability result. First we define the relevant structures.

\begin{definition}   If $\alpha = \alpha(\kappa,\lambda,\mu)$ is the order-type of $Card^-$ in ${\cal
M}_{\kappa \lambda \mu}$, we let
$${\cal N}^2_{\kappa \lambda \mu} \hspace{-.05in} = \hspace{-.05in} ((IS_n)_{n \ge
1},(\alpha_{[n]})_{n \ge 0}, (\alpha^{[n]})_{n \ge 0}; Eq^1, Prod^1, (Proj_n^1)_{n \ge 1}, (<_n)_{n
\ge 0}, (<^n)_{n \ge 0})$$ be the structure whose sorts are viewed as being pairwise disjoint (and
all but finitely many $\alpha_{[n]}$ are empty), and $<_n, <^n$ are the usual (well-) orderings on
$\alpha_{[n]}, \alpha^{[n]}$. The superscript $2$ indicates that ${\cal N}^2_{\kappa \lambda \mu}$ is
viewed as a second order structure in a very strong sense. This means that the language used to
describe it, as well as including first order variables corresponding to each sort, also contains,
for each tuple of sorts, variables ranging over relations whose $i$th entry lies in the $i$th sort of
the tuple for each $i$. (Alternatively we can introduce sorts corresponding to each such tuple,
adjoin all the natural relations, and work in first order logic).   \label{5.4}    \end{definition}

In one direction the interpretability is `explicit'.

\begin{theorem} If $cf(\kappa) > 2^{\aleph_0}$ then ${\cal N}^2_{\kappa \lambda \mu}$ is interpretable
in ${\cal M}^*_{\kappa \lambda \mu}$.           \label{5.5}  \end{theorem}

\noindent{\bf Proof}\hspace{.1in} The main point is to show how each $\alpha_{[n]}, \alpha^{[n]}$ may 
be represented in ${\cal M}^*_{\kappa \lambda \mu}$. Then we sketch how second order variables as
described above are `simulated' within the first order language of ${\cal M}^*_{\kappa \lambda \mu}$.

First we represent non-empty subsets of $Card^-$ of cardinality $< \Omega$ using members of $F_2$, the
idea being that $h \in F_2$ represents its range. Since $2^{\aleph_0} < cf(\kappa) \le \kappa <
\lambda$, all such subsets of $Card^-$ can be represented. The following can then be expressed:

the set encoded by $h_1$ is a subset of that encoded by $h_2$:

\hspace{.2in} $incl(h_1,h_2): (\forall t_1)(\exists t_2)(h_1(t_1) = h_2(t_2))$,

\noindent (where as $App_n$ is just `application' we write $h_1(t_1)$ instead of $App_2(h_1,t_1)$
etc).

$h_1$ and $h_2$ encode the same set: $incl(h_1,h_2) \wedge incl(h_2,h_1)$,

$h$ encodes a final segment of $\alpha$:

\hspace{.2in} $final(h): (\forall t_1)(\forall \beta)(\exists t_2)(h(t_1) \le \beta \rightarrow h(t_2) =
\beta)$,

$\alpha$ is divisible by $\Omega$ ($\equiv \alpha_{[0]} = 0$),

\hspace{.2in} $div(\alpha, \Omega): (\forall h) \neg final(h)$.

$\beta \in \alpha$ is divisible by $\Omega, \Omega^{k+1}$:

\hspace{.2in} $div(\beta, \Omega): (\forall h)(\forall \gamma < \beta)(\exists \delta)(\gamma \le
\delta < \beta \wedge (\forall t)(h(t) \neq \delta))$,

\hspace{.2in} $div(\beta, \Omega^{k+1}): div(\beta, \Omega^k) \wedge (\forall h)(\forall \gamma <
\beta)(\exists \delta)(\gamma \le \delta < \beta \wedge div(\delta,\Omega^k) \wedge
\hspace{3in} (\forall t)(h(t) \neq \delta))$.

$\alpha$ is divisible by $\Omega^{k+1}$:

\hspace{.2in} $div(\alpha, \Omega^{k+1}): div(\alpha, \Omega^k) \wedge (\forall h)(\forall
\beta < \alpha)(\exists \gamma)(\beta \le \gamma < \alpha \wedge div(\gamma, \Omega^k)
\wedge \hspace{3in} (\forall t)(h(t) \neq \gamma))$.

$\alpha_{[0]}$ is now represented by $h$ such that 
$$final(h) \wedge (\forall h')(incl(h,h') \wedge final(h') \rightarrow h = h'),$$ 
if such exists (and otherwise is $0$). Similarly $\alpha_{[k]}$ is represented by $h$ such that
$(\forall t)(div(h(t),\Omega^k)) \wedge (\forall h')((\forall t)(div(h'(t), \Omega^k)) \wedge
incl(h,h') \rightarrow h = h')$ if such exists (and otherwise is $0$).

To encode facts about cofinalities we quantify over non-empty binary relations on $Card^-$ of
cardinality $< \Omega$ using pairs $(h_1,h_2)$ in $F_2$. Observe that if $\emptyset \neq R \subseteq
(Card^-)^2, |R| < \Omega$, then for some $h_1, h_2 \in F_2, R = \{ (h_1(t), h_2(t)): t \in IS_2 \}$.
We can describe when $R$ is an order-isomorphism thus:

\hspace{.2in} $iso(h_1,h_2): (\forall t_1)(\forall t_2)(h_1(t_1) \le h_1(t_2) \leftrightarrow h_2(t_1)
\le h_2(t_2))$.

The set coded by $h$ is then cofinal in $\alpha$ if

\hspace{.2in} $cofinal(h): (\forall \beta)(\exists t)(\beta \le h(t))$,

and $h$ codes the cofinality of $\alpha$, which is $< \Omega$ if

\hspace{.2in} $cofinal(h) \wedge (\forall h')(cofinal(h') \rightarrow (\exists h_1)(\exists
h_2)(iso(h_1,h_2) \wedge incl(h,h_1)$

\hspace{.5in}$\wedge \; incl(h_1, h) \wedge incl(h_2, h'))$.

We may express $cf(\alpha) \ge \Omega$ by $(\forall h) \neg cofinal(h)$, and in a similar way for
each $\beta \in \alpha$ we may express `$h$ codes a cofinal subset of $\beta$' and `$h$ codes
$cf(\beta)$'. From this it should be clear that each $\alpha_{[n]}, \alpha^{[n]}$ can be represented
(though presumably not uniformly).

Finally we show how to represent non-empty $n-$ary relations of cardinality $< \Omega$ on the sorts of
the original structure according to sort provisos of the kind described above. Since the sorts of
${\cal N}^2_{\kappa \lambda \mu}$ all have cardinality $< \Omega$, this translates into {\em full}
second order logic in this structure.

Without loss of generality consider a tuple of sorts in ${\cal M}^*_{\kappa \lambda \mu}$ of the form
$(IS_{i_1}, \ldots, IS_{i_m}, Card^-, \ldots,Card^-)$ (where $Card^-$ occurs $n$ times). We represent
corresponding non-empty relations of cardinality $< \Omega$ by $m+n-$tuples of the form $(h_1, \ldots
,h_m, h_1', \ldots , h_n')$ where $h_j \in F_{i_j+2}$ satisfies $fun_{2,i_j}(h_j)$ and $h_j' \in
F_2$. Such an $m+n-$tuple represents 
$$B = B({\overline h}) = \{ (H_1(t), \ldots , H_m(t), h'_1(t), \ldots , h'_n(t)):  t \in IS_2\}$$
where $H_j$ is the function from $IS_2$ to $IS_{i_j}$ determined by $h_j$.

Clearly $B({\overline h})$ {\em is} a non-empty relation of the required kind of cardinality $<
\Omega$, and conversely every such relation can be written as $B({\overline h})$ for some
$m+n$-tuple ${\overline h}$. 

As in the proof of Lemma \ref{5.1}, ${\overline t} \in B({\overline h})$ can be expressed in ${\cal
M}^*_{\kappa \lambda \mu}$.   

If $R$ is a non-empty $n$-ary relation on ${\cal N}^2_{\kappa \lambda \mu}$ with specified sorts, then
as each individual sort is definable as indicated above, $R$ may be represented by a corresponding
$n$-ary relation of the kind just discussed, in ${\cal M}^*_{\kappa \lambda \mu}$. $\Box$

\vspace{.1in}

In the other direction we have a weaker notion than `semi-interpretability', which is nevertheless
sufficient to transfer elementary equivalence. The weakening just consists in having a whole family
of representatives of a tuple rather than a single one. Let us say that for $k \in {\Bbb N}$ a {\em
$k$-representation} of a tuple $({\overline t}, {\overline \beta}, {\overline h})$ in ${\cal
M}^*_{\kappa \lambda \mu}$ where each $t_i$ lies in some $IS_j, \beta_i \in Card^-$, and each $h_i$
lies in some $F_j$, is any tuple of the form $(A, <^*, {\overline g}, {\overline t},  {\overline b},
{\overline H})$ where 

$<^*$ is a well-ordering of $IS_2$,

$A \subseteq IS_2, b_i \in A$,

if $h_i \in F_j$ then $H_i: IS_j \rightarrow A$,

${\overline g} = (g_0, \ldots, g_{k-1}, g^0, \ldots , g^{k-1})$ where $g_i, g^i: A \cup \{ \infty \}
\rightarrow IS_2$,

for some order-preserving 1--1 map $\theta : A \rightarrow Card^-, \theta (b_i) = \beta_i$,
$(\forall t \in IS_j) \theta (H_i(t)) = h_i(t)$, and for every $a \in A \cup \{ \infty \}$ the
order-type of $\{ t \in IS_2: t <^* g_j(a) \}$ equals $\gamma(\theta(a),\theta(A))_{[j]}$ and the
order-type of $\{ t \in IS_2: t <^* g^j(a) \}$ equals $\gamma(\theta(a),\theta(A))^{[j]}$ (where we
take $\theta(\infty) = \lambda \; (> \beta$ for all $\beta \in Card^-$)).

We remark that all entries in this tuple except for the $g_i, g^i$ lie in ${\cal N}^2_{\kappa \lambda
\mu}$:
$$<^* \subseteq IS_2^2, A \subseteq IS_2, t_i \in IS_j, b_i \in IS_2, H_i \subseteq IS_j
\times IS_2.$$
Moreover $g_i|A, g^i|A \subseteq IS_2^2$ lie in ${\cal N}^2_{\kappa \lambda \mu}$ so by making an easy
modification to their `official' definition, so do $g_i, g^i$. But $\theta$ does not (which is why it
does not form part of the representation).

\begin{lemma} Any $({\overline t}, {\overline \beta}, {\overline h})$ has a $k$-representation. 
\label{5.6}   \end{lemma}

\noindent {\bf Proof} \hspace{.1in} Let $A'$ be the union of the set of entries of $\overline
\beta$ and the ranges of the $h_i, A' \subseteq Card^-$. Then $|A'| \le 2^{\aleph_0}$. Also
$\gamma(a,A')_{[j]}, \gamma(a,A')^{[j]}$ each has order-type at most that of $A'$. We choose $A
\subseteq IS_2$ of cardinality $|A'|$, a bijection $\theta: A \rightarrow A'$, and a well-ordering
$<^*$ of $IS_2$ extending $\theta^{-1}(<)$. For $a \in A$ let $g_j(a)$ equal the
$\gamma(\theta (a), \theta (A))_{[j]}$th element of $IS_2$ under $<^*$, $g_j(\infty) = \lambda$, and
similarly for $g^j(a)$. Let $b_i = \theta^{-1}(\beta_i)$ and $H_i(t) = \theta^{-1}h_i(t)$ for each
$t$.     $\Box$

\begin{lemma} In the language of ${\cal N}^2_{\kappa \lambda \mu}$ for each $k$ there is a
formula $\varphi_k$ such that ${\cal N}^2_{\kappa \lambda \mu} \models \varphi_k[A, <^*,{\overline
g},b]$ if and only if

$A \subseteq IS_2, <^*$ is a well-ordering of $IS_2, b \in IS_2$,

${\overline g} = (g_0, \ldots , g_{k-1}, g^0, \ldots , g^{k-1})$ where the $g_i, g^i$ are
functions from $A \cup \{ \infty \}$ into $IS_2$,

and if $A$ and $IS_2$ are enumerated in $<^*$-increasing order as $\{ a_\beta: \beta < \beta_0 \},
\{ b_\gamma: \gamma < \gamma_0 \}$, and $b = b_{\gamma_1}$, and for each $\beta < \beta_0,
\alpha_\beta$ is an ordinal for which $b_{(\alpha_\beta)_{[l]}} = g_l(a_\beta),
b_{(\alpha_\beta)^{[l]}} = g^l(a_\beta)$ for $l < k$, then $\sum_{\beta <
\beta_0}\alpha_\beta \sim_k \gamma_1$.       \label{5.7}   \end{lemma}

\noindent{\bf Proof}\hspace{.1in} The proof of this is obtained by formalizing a transfinite induction
similar to that used in the proof of Lemma \ref{5.3}(iii).    $\Box$

\begin{lemma} For each tuple of sorts and each $k$, 

(i) there is a formula $rep_k$ of the language of ${\cal N}^2_{\kappa \lambda \mu}$ which holds in 
${\cal N}^2_{\kappa \lambda \mu}$ for a tuple having the right sequence of sorts if and only if it is
a $k$-representation (of some tuple),

(ii) there is a formula $isorep_k$ of the language of ${\cal N}^2_{\kappa \lambda \mu}$ which holds
in  ${\cal N}^2_{\kappa \lambda \mu}$ for a pair of tuples each having the right sequence of sorts if
and only if there is some tuple of ${\cal M}^*_{\kappa \lambda \mu}$ of which they are both
$k$-representations.       \label{5.8}   \end{lemma}

\noindent{\bf Proof}\hspace{.1in} (i) To tell whether a tuple is a $k$-representation we first verify lines $1$
to $4$ of the definition, which can all be expressed in the language of ${\cal N}^2_{\kappa \lambda
\mu}$ (where we have second order logic). If they hold then the main point is to check whether
$\theta$ can be defined to give the correct $\gamma$-values. For this we appeal to the previous
lemma, and we also need to refer to the sorts $\alpha_{[j]}, \alpha^{[j]}$ for $j < k$ to ensure that
the right $\gamma(\theta(\infty), \theta(A))_{[j]}, \gamma(\theta(\infty), \theta(A))^{[j]}$ values
can be achieved. Then we may define ${\overline \beta}, {\overline h}$ by $\beta_i = \theta(b_i),
h_i(t) = \theta(H_i(t))$.

(ii) Similar remarks apply except that we should now work with `minimal' $A$, that is, those
which are equal to the union of the $\{ b_i \}$ and $range(H_i)$.   $\Box$ 

\begin{theorem} For every (first order) formula $\varphi(x_0, \ldots , x_{n-1})$ of the
language of ${\cal M}^*_{\kappa \lambda \mu}$ there is an effectively determined integer $k$ and
(second order) formula $\psi(y_0, \ldots, y_{2k+n+1})$ of the language of ${\cal N}^2_{\kappa \lambda
\mu}$ such that for all $\kappa, \lambda,\mu$ with $cf(\kappa) > 2^{\aleph_0}$, for every $a_0,\ldots,
a_{n-1}$ in ${\cal M}^*_{\kappa \lambda \mu}$ (having the correct sorts) and every 
$k$-representation $c_0, \ldots , c_{2k+n+1}$ of $\overline a$ in ${\cal N}^2_{\kappa \lambda \mu}$,
$${\cal M}^*_{\kappa \lambda \mu} \models \varphi[{\overline a}] \Leftrightarrow {\cal N}^2_{\kappa
\lambda \mu} \models \psi[{\overline c}].$$               
                               \label{5.9}      \end{theorem}

\noindent{\bf Proof}\hspace{.1in} We construct $\psi$ by induction. The $k$ is just the `quantifier 
depth' of $\varphi$ (for quantifications over $F_j$), as emerges from what follows.

First consider the case of atomic formulae, where we take $k = 0$. If $\varphi(x_0)$ is $Eq^1(x_0)$
we let $\psi(y_0)$ also be $Eq^1(y_2)$, (since sort $IS_2$ is the same in the two structures).
Similarly for $Prod^1(x_0), Proj_n^1(x_0,x_1)$, and $x_0 = x_1$ where $x_0, x_1$ lie in the same
$IS_j$. If $\varphi(x_0,x_1)$ is $x_0 = x_1$ or $x_0 < x_1$ where $x_0, x_1 \in Card^-$ we let
$\psi(y_0, y_1, y_2, y_3)$ be $y_2 = y_3$ or $(y_2,y_3) \in y_1$ respectively. Consider $<$ for
instance, and let $(A, <^*, b_0, b_1)$, a $k$-representation of $(\beta_0, \beta_1 )$, and $\theta$,
be given by the definition of what this means. Then $\theta(b_i) = \beta_i$, so $b_0 <^* b_1
\Leftrightarrow \beta_0 < \beta_1$ and ${\cal M}^*_{\kappa \lambda \mu} \models
\varphi[\beta_0,\beta_1] \Leftrightarrow \beta_0 < \beta_1 \Leftrightarrow b_0 <^* b_1
\Leftrightarrow {\cal N}^2_{\kappa \lambda \mu} \models (b_0,b_1) \in <^* \Leftrightarrow  {\cal
N}^2_{\kappa \lambda \mu} \models \psi[A,<^*,b_0,b_1]$. For the remainder it suffices to consider
$\varphi(x_0,x_1,x_2) \equiv x_1 = App_n(x_2,x_0)$ (since the other atomic formulae may be written
in terms of this and the ones above). Here we let $\psi(y_0,y_1,y_2,y_3,y_4)$ be $(y_2,y_3) \in y_4$.

For the induction step the case of negation is immediate (we take the same $k$ and the negation of
the corresponding formula). For conjunction suppose that $\varphi$ is $\varphi_1(x_0, ..
,x_{l-1},x_l, .. , x_{m-1}) \wedge \varphi_2(x_0, \ldots ,x_{l-1},x_m, \ldots , x_{n-1})$,
where the $x_i$ are distinct variables , and that $k_1, \psi_1(y_0', .. , y_{2k_1+1}', x_0',
. , x_{m-1}')$ corresponding to $\varphi_1$ and $k_2, \psi_2(y_0'', \ldots , y_{2k_2+1}'',
x_0'', \ldots , x_{l-1}'', x_m'', \ldots , x_{n-1}'')$ corresponding to $\varphi_2$ have been
chosen. Let $k = max(k_1,k_2)$ and $\psi(y_0, \ldots , y_{2k+1}, z_0, \ldots , z_{l-1}, z_l, \ldots,$

\noindent $z_{m-1}, z_m, \ldots , z_{n-1})$ be the formula
$$\psi_1(y_0, \ldots , y_{2k_1+1}, z_0, \ldots , z_{m-1}) \wedge \psi_2(y_0, \ldots ,
y_{2k_2+1}, z_0, \ldots , z_{l-1}, z_m, \ldots , z_{n-1}).$$
If $a_0, \ldots , a_{n-1}$ in ${\cal M}^*_{\kappa \lambda \mu}$ have the correct sorts, and
$(A, <^*, g_0, \ldots , g_{k-1}, g^0, \ldots ,$

\noindent $g^{k-1}, b_0, \ldots , b_{n-1})$ is a $k$-representation in ${\cal
N}^2_{\kappa \lambda \mu}$, then $(A, <^*, g_0, \ldots, g_{k_1-1}, g^0,$

\noindent $\ldots , g^{k_1-1}, b_0, \ldots , b_{m-1}), (A, <^*, g_0, \ldots ,
g_{k_2-1}, g^0, \ldots, g^{k_2-1}, b_0, \ldots , b_{l-1}, b_m, \ldots ,$

\noindent $b_{n-1})$ are $k_1$-,$k_2$-representations of $(a_0, \ldots , a_{m-1}), (a_0,
\ldots , a_{l-1}, a_m, \ldots , a_{n-1})$

\noindent respectively, and so the result goes through with this $k$.

Now consider the existential quantifier. Suppose $\varphi(x_0, \ldots , x_{n-1})$ is $(\exists x_n)$

\noindent $\varphi'(x_0, \ldots , x_{n-1}, x_n)$, and that $k', \psi'$ corresponding to
$\varphi'$ have been chosen.

\noindent{\bf Case 1}: $x_n \in IS_j$. Let $k = k'$ and $\psi(y_0, \ldots , y_{2k+n+1})$ be $(\exists
y_{2k+n+2}) \psi'(y_0, \ldots ,$

\noindent $y_{2k+n+2})$ (where $y_{2k+n+2} \in IS_j$ too).

Suppose $a_0, \ldots , a_{n-1} \in {\cal M}^*_{\kappa \lambda \mu}$ have the correct sorts, and $(A,
<^*, {\overline g}, {\overline c})$ is a $k$-representation of ${\overline a}$. Then for any $a_n
\in IS_j, (A, <^*, {\overline g}, {\overline c}, a_n)$ is a $k$-representation of $(a_0, \ldots ,
a_n)$. Hence
 
\vspace{.1in}
${\cal M}^*_{\kappa \lambda \mu} \models \varphi[a_0, \ldots , a_{n-1}]$
\vspace{-.1in}
\begin{eqnarray*} & \Leftrightarrow & \mbox{for some } a_n \in IS_j, {\cal M}^*_{\kappa \lambda \mu}
\models \varphi'[a_0, \ldots , a_n] \\   
& \Leftrightarrow & \mbox{ for some } c_n \in IS_j, {\cal N}^2_{\kappa \lambda \mu} \models
\psi'[A, <^*, {\overline g}, c_0, \ldots , c_n]\\ 
& \Leftrightarrow & {\cal N}^2_{\kappa \lambda \mu} \models
\psi[A, <^*, {\overline g}, {\overline c}].
\end{eqnarray*}

\noindent{\bf Case 2}: $x_n \in Card^-$. Any existential quantifiers over $Card^-$ may be eliminated
in favour of quantifiers over $IS_2$ and $F_2$, since $(\exists x_n \in Card^-) \varphi'(x_0, \ldots
, x_n) \Leftrightarrow$

\noindent $(\exists h \in F_2)(\exists t \in IS_2) \varphi'(x_0, \ldots , x_{n-1}, h(t))$.

\noindent{\bf Case 3}: $x_n \in F_j$. Let $k = k'+1$ and $\psi(y_0,\ldots,y_{2k+n+1})$ be the formula
$$\begin{array}{c}
(\exists z_0) \ldots (\exists z_{2k+n+1})(\exists z)(rep_k(z_0, \ldots , z_{2k+n+1},z) \; \wedge \\
isorep_k(y_0, \ldots , y_{2k+n+1},z_0, \ldots , z_{2k+n+1}) \wedge \psi''(z_0, \ldots ,
z_{2k+n+1},z)),
     \end{array} $$
where $rep_k, isorep_k$ are the appropriate instances of the formulae provided by Lemma \ref{5.8}
(that is for the correct sequence of sorts), and $\psi''(z_0, \ldots , z_{2k+n+1},$

\noindent $z)$ is $\psi'(z'_0, \ldots , z'_{2k+n-1},z)$ where $(z'_0, \ldots ,
z'_{2k+n-1},z)$ is obtained from $(z_0, \ldots ,$

\noindent $z_{2k+n+1},z)$ by deleting the two variables corresponding to $g_k$ and $g^k$. For ease
assume the variables in $\bigcup IS_i$ come first, then those in $Card^-$, then those in $\bigcup
F_i$.

Let $(A,<^*,{\overline g}, {\overline t}, {\overline b}, {\overline H})$ be a $k$-representation of
$({\overline t}, {\overline \beta}, {\overline h})$. Then
\vspace{.1in} 

\noindent ${\cal M}^*_{\kappa \lambda \mu}
\models \varphi[{\overline t}, {\overline \beta}, {\overline h}]$
\vspace{-.1in} 
\begin{eqnarray*}
& \Leftrightarrow & \mbox{ for some } h \in F_j, {\cal M}^*_{\kappa \lambda \mu} \models
\varphi'[{\overline t}, {\overline \beta}, {\overline h}, h] \\    
& \Leftrightarrow & \mbox{for some } A', <', {\overline {g'}}, {\overline t}, {\overline {b'}},
{\overline {H'}}, H', \mbox{ where } H': IS_j \rightarrow A', \\  
&  & {\cal N}^2_{\kappa \lambda \mu} \models rep_k[A', <' ,{\overline {g'}}, {\overline t},
{\overline {b'}}, {\overline {H'}}, H'] \wedge \\  
&  & isorep_k[A', <' , {\overline {g'}}, {\overline t}, {\overline {b'}}, {\overline {H'}}, A, <^*,
{\overline g}, {\overline t}, {\overline b}, {\overline H}] \wedge \psi''[A', <', {\overline {g'}},
{\overline t}, {\overline {b'}}, {\overline {H'}}, H'] \\    
& \Leftrightarrow & {\cal N}^2_{\kappa \lambda \mu} \models \psi[A, <^*, {\overline g},  {\overline
t}, {\overline b}, {\overline H}].  \end{eqnarray*} 

The first and last steps are immediate. It is the intermediate equivalence which we have to justify.

Suppose then that ${\cal M}^*_{\kappa \lambda \mu} \models \varphi'[{\overline t}, {\overline \beta},
{\overline h}, h]$, and let $(A', <' , {\overline {g'}}, {\overline t}, {\overline {b'}}, {\overline
{H'}}, H')$ be a $k$-representation of $({\overline t}, {\overline \beta}, {\overline h}, h)$ (which
exists by Lemma \ref{5.6}). We get a corresponding $k'$-representation by omitting $g_k',
(g^k)'$, so by the induction hypothesis, ${\cal N}^2_{\kappa \lambda \mu} \models \psi''[A', <',
{\overline {g'}},  {\overline t}, {\overline {b'}}, {\overline {H'}}, H']$. Also ${\cal N}^2_{\kappa
\lambda \mu} \models rep_k[A', <', {\overline g'},  {\overline t}, {\overline {b'}}, {\overline
{H'}}, H'] \wedge$

\noindent $isorep_k[A', <', {\overline {g'}},  {\overline t}, {\overline {b'}}, {\overline {H'}}, A,
<^*, {\overline g},  {\overline t}, {\overline b}, {\overline H}]$.

Conversely if ${\cal N}^2_{\kappa \lambda \mu} \models \psi[A, <^*, {\overline g},  {\overline
t}, {\overline b}, {\overline H}]$ there are $A', <', {\overline {g'}},  {\overline t}, 
{\overline {b'}}, {\overline {H'}}, H'$ which form a $k$-representation of ${\overline t}, {\overline
\beta}, {\overline h}, h$ for some $h$, and such that ${\cal N}^2_{\kappa \lambda \mu} \models
rep_k[A', <'~,$

\noindent ${\overline g'},  {\overline t}, {\overline {b'}}, {\overline {H'}}, H'] \wedge \;
isorep_k[A', <', {\overline {g'}},  {\overline t}, {\overline {b'}}, {\overline {H'}}, A, <^*,
{\overline g},  {\overline t}, {\overline b}, {\overline H}]$. By Lemma \ref{5.3}(iv) there is
some $H: IS_j \rightarrow A$ such that $(A, <^*, {\overline g},  {\overline t}, {\overline b},
{\overline H}, H)$ is a $k'$-representation of  $({\overline t}, {\overline \beta}, {\overline h},
h)$. By the induction hypothesis again,  ${\cal M}^*_{\kappa \lambda \mu} \models \varphi'[{\overline
t}, {\overline \beta}, {\overline h},h]$ as required.    $\Box$

\begin{corollary} If $cf(\kappa_1), cf(\kappa_2) > 2^{\aleph_0}$ then the following are equivalent:

(i) ${\cal M}_{\kappa_1 \lambda_1 \mu_1} \equiv {\cal M}_{\kappa_2 \lambda_2 \mu_2}$, 

(ii) ${\cal N}^2_{\kappa_1 \lambda_1 \mu_1} \equiv {\cal N}^2_{\kappa_2 \lambda_2 \mu_2}$,

(iii) $S_{\lambda_1}(\mu_1)/S_{\kappa_1}(\mu_1) \equiv S_{\lambda_2}(\mu_2)/S_{\kappa_2}(\mu_2)$. 
                                      \label{5.10}      \end{corollary}

\noindent{\bf Proof}\hspace{.1in} This follows from Theorems \ref{5.2}, \ref{5.5}, and \ref{5.9}.
 $\Box$ 

\vspace{.1in}

So in a certain sense, for cofinalities above $2^{\aleph_0}$, only a rather modest amount of
information about the cardinals $\kappa, \lambda$, and $\mu$ is needed to distinguish the quotient
groups, and in particular, whenever $\alpha(\kappa_1,\lambda_1,\mu_1) \sim
\alpha(\kappa_2,\lambda_2,\mu_2)$ they are elementarily equivalent.

\section{The case $cf(\kappa) \le 2^{\aleph_0}$ and conclusions}

In this section we begin by treating the rather more complicated case in which $cf(\kappa) \le
2^{\aleph_0}$, and then summarize the conclusions in all cases. The first remark is that there is a
first order sentence of the language of ${\cal M}_{\kappa \lambda \mu}$ which distinguishes this case,
namely
$$(\exists h_1, h_2 \in F_2)(h_1 \neq h_2 \wedge (\forall t \in IS_2)(App_2(h_1,t) =
App_2(h_2,t))).$$
So from now on we assume that $cf(\kappa) \le 2^{\aleph_0}$.

We now describe the modification of ${\cal N}^2_{\kappa \lambda \mu}$ appropriate in this case, which
varies slightly according as $\kappa \le 2^{\aleph_0}$ or not, and $\kappa = \aleph_0$ or not (cases
which we shall see below can be distinguished by formulae of the language of group theory). Let
$\alpha^*$ be the least ordinal $> 0$ such that $(\exists \gamma)(\beta = \gamma + \alpha^*)$ where
$\kappa = \aleph_\beta$. The definition of ${\cal N}^2_{\kappa \lambda \mu}$ is modified to include as
additional sorts $cf(\kappa)$, and $\alpha^*_{[n]}, \alpha^{*[n]}$ for $n \ge 0$. Since $\alpha^*$ is
by definition additively indecomposable, only at most one $\alpha^*_{[n]}$ can be non-zero, so the
representation is somewhat redundant, and we have just $\omega + 1$ possible cases. We also include
(distinct) individual constants $c_0, c_\kappa \in IS_2$ in the structure. These may be chosen
arbitrarily or, better, as {\em definable} elements (to ensure that the interpretation is without
parameters). 

\begin{definition}   If $\alpha = \alpha(\kappa,\lambda,\mu)$ is the order-type of $Card^-$ in ${\cal
M}_{\kappa \lambda \mu}$ and $\alpha^*$ is the least ordinal $> 0$ such that $(\exists \gamma)(\beta =
\gamma + \alpha^*)$ where $\kappa = \aleph_\beta$, we let
$$\begin{array}{c}
{\cal N}^2_{\kappa \lambda \mu} = ((IS_n)_{n \ge 1}, (\alpha_{[n]})_{n \ge 0}, (\alpha^{[n]})_{n \ge
0},(\alpha^*_{[n]})_{n \ge 0}, (\alpha^{*[n]})_{n \ge 0}; Eq^1, Prod^1, \\
 (Proj_n^1)_{n \ge 1}, (<_n)_{n \ge 0}, (<^n)_{n \ge 0},(<_n^*)_{n \ge 0},(<^{*n})_{n \ge 0}, c_0,
c_\kappa, kap, fin)
     \end{array} $$
be the structure whose sorts are viewed as being pairwise disjoint (and all but finitely many
$\alpha_{[n]}$ and all but at most one $\alpha^*_{[n]}$ are empty), and $<_n, <^n, <_n^*, <^{*n}$ are
the usual (well-) orderings on $\alpha_{[n]}, \alpha^{[n]}, \alpha^*_{[n]}, \alpha^{*[n]}$. As in
Definition \ref{5.4} the superscript $2$ indicates that ${\cal N}^2_{\kappa \lambda \mu}$ is a second
order structure, and the same restrictions are made as before on the second order variables which are
allowed (where the new sorts are now allowed as entries in the tuples of sorts), except that since we
no longer know for sure that $\lambda > 2^{\aleph_0}$, we have to restrict to quantification over
relations of cardinality $< \lambda$. The constants $c_0$ and $c_\kappa$ are distinct elements of
$IS_2$, and $kap$ and $fin$ are unary relations on $IS_2$, $kap$ picking out a subset of $IS_2$ of
cardinality $\kappa$ and $fin$ the set of isomorphism types of finite sets, which are only included
if $\kappa \le 2^{\aleph_0}$, $\kappa = \aleph_0$ respectively.   \label{6.1}   \end{definition}

The case $\kappa \le 2^{\aleph_0}$ has to be treated separately because it is precisely here that
$Sum_n$ (summation of $h \in F_n$) cannot be identified with supremum. As we saw above, subsets of
$IS_2$ can be represented in ${\cal M}_{\kappa \lambda \mu}$, and by various tricks (which we do not
go into, but which are similar to ones described below for other purposes) one can express the
property of having cardinality $\kappa$. In general there will be no {\em definable} such set
however, so the interpretation of ${\cal N}^2_{\kappa \lambda \mu}$ in ${\cal M}_{\kappa \lambda \mu}$
in this case requires a parameter. If also $\kappa = \aleph_0$, we include a predicate $fin$ picking
out out the members of $IS_2$ corresponding to isomorphism types of finite sets. (This predicate {\em
is} definable in ${\cal M}_{\kappa \lambda \mu}$.) 

The fact that $kap$ is not definable does not affect our main results however. We shall show that
(in the relevant case), $kap$ can be interpreted, and that we can express when the representations
of members of our structure using two possible interpretations of $kap$ represent the same object.  

One main difference in this section is that we can no longer work with ${\cal M}^*_{\kappa \lambda
\mu}$. Instead we refine the methods of section 5 to show how the second order logic just mentioned
can be represented in ${\cal M}_{\kappa \lambda \mu}$. We recall that in Lemma \ref{5.1}(ii) we
saw how to say that two members of $F_n$ or $F_{m+n}$ encode the same subsets of $IS_n$ or $IS_m
\times IS_n$. In fact if $cf(\kappa) > 2^{\aleph_0}$ they encode the same set if and only if they
are equal. But this is not true if $cf(\kappa) \le 2^{\aleph_0}$ (as was essentially exploited above
in devising a sentence to characterize this case). Life is easier if we use $h$ which `minimally
encode' sets or relations. All this means is that the cumulative effect of values below $\kappa$ is
negligible, in other words $\sum \{h(t): h(t) < \kappa \} < \kappa$, but we have to see how this can
be formally expressed. 

For $h_1, h_2 \in F_n$ we write $restr_n(h_1,h_2)$ for $\sum \{ h_1(t)-h_2(t): t \in IS_n \} <
\kappa$. In $S_\lambda(\mu)/S_\kappa(\mu)$ this corresponds to a tuple representing $h_1$ being
conjugate to a restriction of a tuple representing $h_2$ (expressed in section 4 by a corresponding
formula $restr_n$), and so by Theorem \ref{2.6} is first order expressible in the language of ${\cal
M}_{\kappa \lambda \mu}$. Saying that $h$ minimally encodes a set (or relation) then is expressed by
$min(h)$:
$$(\forall h' \in F_n)((\forall t \in IS_n)(App_n(h,t) = \kappa \leftrightarrow App_n(h',t) =
\kappa) \rightarrow restr_n(h,h')).$$

Now we show how to capture the behaviour of cardinals below $\kappa$ in ${\cal M}_{\kappa \lambda
\mu}$. Let us write $Card_{< \kappa}$ for $\{\nu \in Card: \nu < \kappa \}$. We can only hope to
capture the `tail' of $Card_{< \kappa}$. We encode (the tail of) a subset $X$ of $Card_{<
\kappa}$ by any $k \in F_2$ having $X$ as range. (Of course subsets of $Card_{< \kappa}$ of
cardinality $< \min(\Omega, \lambda)$ can be so encoded.) We can express `$k$ encodes some set' by
$(\forall t \in IS_2)(App_2(k,t) = 0)$, and we say that such a $k$ is {\em almost zero}. In the sense
of the previous paragraph $k$ encodes the empty subset of $IS_2$. As we wish to exclude $0$ (that is,
any $k$ such that $\sum \{k(t): t \in IS_2\} < \kappa)$ we identify $0$ as any $k \in F_2$ which {\em
minimally} encodes the empty set.

In order to express when two almost zero members of $F_2$ encode the same subset of $Card_{< \kappa}$
it is easier to pass to those which are `almost 1--1', meaning that 
$$(\exists \nu < \kappa)(\forall t_1, t_2 \in IS_2)(k(t_1) = k(t_2) \ge \nu \rightarrow t_1 = t_2).$$
This requires a further technical trick.

Now if $h \in F_2$ minimally encodes a subset $X$ of $IS_2$, and $k_1, k_2$ are almost zero, we can
express `$k_1$ encodes the restriction of $k_2$ to $X$' by the formula
$$restr_2(k_1,k_2) \wedge restr_2(k_1,h) \wedge (\forall k' \in F_2)(restr_2(k',k_2) \wedge
restr_2(k',h) \rightarrow restr_2(k',k_1)).$$
If $f \in F_4$ and $one$-$onefun_{2,2}(f)$ we can express `the function $F$ coded by $f$ carries
$k_1$ to $k_2$' (meaning that $F$ carries $\{t: k_1(t) > 0\}$ to $\{t: k_2(t) > 0\}$, and for each
$t$ with $k_1(t) > 0, k_2(F(t)) = k_1(t)$), via $S_\lambda(\mu)/S_\kappa(\mu)$ and Theorem \ref{2.6}
as follows:

`a tuple representing $f$ has a restriction which projects to a conjugate of a tuple representing
$k_1$ on co-ordinates 1,2, and to a conjugate of a tuple representing $k_2$ on co-ordinates 5,6'.

Using this we can now express `$k$ is almost 1--1' thus:

$one$-$one(k): k$ is almost zero and $\forall k_1 \forall k_2 \forall f (k_1, k_2$ non-zero 
restrictions of $k$ to disjoint subsets of $IS_2 \wedge one$-$onefun(f) \rightarrow \neg (f$
carries $k_1$ to $k_2))$.

For if $k$ is not almost 1--1 there are cofinally many $\nu < \kappa$ such that $|k^{-1}(\nu)| \ge
2$ and we can find non-zero restrictions of $k$ to disjoint subsets of $IS_2$ and a permutation
taking one to the other.

The point of doing this is that we can now express `almost zero $k_1$ and $k_2$ code the same (tail
of a) subset of $Card_{< \kappa}$', and compare order-types of such subsets. For $k_1$ and $k_2$
code the same subset of $Card_{< \kappa}$ if and only if one can be carried to the other by a 1--1
function from a subset of $IS_2$ to $IS_2$.

We can now express $cf(\kappa) \le 2^{\aleph_0} \wedge \kappa$ is a successor by
$$cf(\kappa) \le 2^{\aleph_0} \mbox{ (already expressed) } \wedge (\forall k)(\mbox{$one$-$one$}(k)
\rightarrow k = 0),$$
if desired (though it corresponds to the special case $\alpha^* = 1$).

Now suppose that $cf(\kappa) \le 2^{\aleph_0} \wedge \kappa$ is a limit. We wish to represent
$cf(\kappa)$ and each $\alpha^*_{[n]}$ and $\alpha^{*[n]}$ in ${\cal M}_{\kappa \lambda \mu}$. We
represent $cf(\kappa)$ by any $k$ such that

$one$-$one(k) \wedge k \neq 0 \wedge (\forall k')(one$-$one(k') \wedge k' \neq 0 \rightarrow (\exists
g)(g$ a 1--1 map from a subset of $IS_2$ into $IS_2 \wedge (\forall t)(k(t) \le k'(gt)))$.

For this we need to express $(\forall t)(k(t) \le k'(gt)))$ in ${\cal M}_{\kappa \lambda \mu}$, and
we use the same idea as above, going via $S_\lambda(\mu)/S_\kappa(\mu)$, and say that the
projection to co-ordinates 1,2 of a tuple representing $g$ has as a restriction a conjugate of $k_1$.

Now moving towards representing the $\alpha_{[n]}^*$ and $\alpha^{*[n]}$, we find a formula
$subset^*(k_1,k_2)$ which expresses `$k_1,k_2$ are almost 1--1, and the set encoded by $k_1$ is a
subset of the set encoded by $k_2$' thus:
$$\begin{array}{c}
(\exists f)(\exists h)(\mbox{$one$-$one$}function (f) \wedge h \mbox{ codes a subset of } IS_2
\; \wedge    \\
k_1 \mbox{ is the restriction of $k_2g$ to the set encoded by } h).
     \end{array} $$
To represent $\alpha_{[n]}^*$ in ${\cal M}_{\kappa \lambda \mu}$, the main point is to find
inductively a formula $div^*(h, \Omega^n)$ analogous to the $div$ formulae considered
earlier, expressing `$h$ encodes a function from $IS_2$ to $Card_{< \kappa}$ such that for every
$t$, $h(t)$ is divisible by $\Omega^n$'. For the basis case $div^*(h, \Omega^0)$ just says that
$h$ encodes a function from $IS_2$ to $Card_{< \kappa}$, in other words, $h$ is `almost zero'. We
also need similar almost zero functions from $IS_2^2$ to $Card_{< \kappa}$.

Assuming inductively that $div^*(h, \Omega^n)$ has been found, we take for $div^*(h,
\Omega^{n+1})$ the formula
$$\begin{array}{c}
div^*(h, \Omega^n) \wedge (\forall h')[(h' \mbox{ codes a function from } IS_2^2 \mbox{ to }
Card_{< \kappa}) \wedge  \\
(\forall t,t' \in IS_2)(h'(t,t') < h(t)) \rightarrow (\exists h'')(h'' \mbox{ codes a function from
$IS_2$ to }  \\
Card_{< \kappa} \wedge (\forall t,t' \in IS_2)(h'(t,t') < h''(t) < h(t))].
     \end{array} $$
We illustrated how to handle inequalities in this context above, so such a formula exists, and is
clearly as required.

We can therefore represent each $\alpha_{[n]}^*$ in ${\cal M}_{\kappa \lambda \mu}$. Moreover, if
$\alpha_{[n]}^* \neq 0$ for some $n$, $\alpha^* = \Omega^n.\alpha_{[n]}^*$, and all cofinalities
are at once represented (equal to either 0 or $cf(\Omega^n.\alpha_{[n]}^*)$), and if $\alpha_{[n]}^* =
0$ for all $n$, $\alpha^* = \Omega^\omega.\alpha_{\omega}^*$, so the cofinalities are all equal
to $cf(\kappa)$. Thus all the sorts of ${\cal N}_{\kappa \lambda \mu}^2$ are represented. The
method for representing the second order logic on ${\cal N}_{\kappa \lambda \mu}^2$ described
above is as in the proof of Theorem \ref{5.5}. 

Next we show how to handle the case $\kappa \le 2^{\aleph_0}$. Let us say that $h \in F_n$ {\em takes
at most two values} if for some $h': IS_2 \rightarrow Card_{< \kappa}$, $(h')_{{\cal E}_n} = h$ and
$|range \; h'| \le 2$. This notion is captured in ${\cal M}_{\kappa \lambda \mu}$ by the formula
$$(\forall X \subseteq IS_2)(\exists Y \subseteq X)(\mbox{all permutations of $X$ fixing $Y$
setwise also fix } h).$$
Observe that we need the $\forall \exists$ quantification because we can only quantify over subsets
of $IS_2$ of cardinality $< \lambda$, and we have not insisted that $\lambda > 2^{\aleph_0}$.

We can now characterize $\kappa \le 2^{\aleph_0}$ by means of the formula
$$(\exists h)(h \neq 0 \; \wedge \; h \mbox{ is almost zero $\; \wedge \; h$ takes at most two
values}),$$ 
which justifies defining ${\cal N}^2_{\kappa \lambda \mu}$ by the cases $\kappa > 2^{\aleph_0}$ or
$\kappa \le 2^{\aleph_0}$. All the ingredients of this structure have been represented in ${\cal
M}_{\kappa \lambda \mu}$ in the case $cf(\kappa) \le 2^{\aleph_0} < \kappa$, and when $\kappa \le
2^{\aleph_0}$ we interpret $kap$ as a subset of $IS_2$ of cardinality $\kappa$. We remark that in this
case, $|Card_{< \kappa}| \le 2^{\aleph_0}$, and so this is an instance where the $\alpha^*_{[n]}$ and
$\alpha^{[*n]}$ really are mostly redundant, since $\alpha^*_{[0]} = \alpha^*$, and all other
$\alpha^*_{[n]}$ are zero. If $\kappa = \aleph_0$, we also have to represent $fin$, as mentioned
earlier, and this is done as follows. Amplifying the remarks just before Theorem \ref{4.3}, let us say
that an $n$-tuple ${\overline x} \in S_\lambda(\mu)/S_\kappa (\mu)$ is {\em irreducible} if
${\overline x} \neq 1$ and $\forall {\overline y} \forall {\overline z}(disj_n({\overline y},
{\overline z}) \wedge {\overline x} = {\overline y} * {\overline z} \rightarrow ({\overline y} = 1
\vee {\overline z} = 1))$. Then one easily checks that $S_\lambda(\mu)/S_\kappa (\mu) \models (\exists
{\overline x})({\overline x}$ irreducible) $\Leftrightarrow \kappa = \aleph_0$, and so, by Theorem
\ref{2.6}, this can also be expressed in ${\cal M}_{\kappa \lambda \mu}$. Moreover, the same argument
shows that irreducibility too can be expressed in ${\cal M}_{\kappa \lambda \mu}$, and we note that
$t \in fin \Leftrightarrow (\forall h \in F_2)(h$ irreducible $\rightarrow App_2(H,t) = 0)$. For if
$S_\kappa(\mu).{\overline g}$ is irreducible and $App_2(Ch_{\overline g},t) \neq 0$, where $t \in
fin$, then $\langle {\overline g}\rangle$ must have infinitely many orbits of type $t$, so can be
written as a non-trivial product of disjoint elements. On the other hand, if $t \not \in fin$, then
there is $S_\kappa(\mu).{\overline g} \neq 1$ such that $\langle {\overline g}\rangle$ has a single
non-trivial orbit of type $t$.          

This describes the essential steps in the proof of the following theorem.

\begin{theorem} For all $\kappa$, ${\cal N}_{\kappa \lambda \mu}^2$ is interpretable in ${\cal
M}_{\kappa \lambda \mu}$. If $\kappa > 2^{\aleph_0}$ the interpretation is without parameters,
and if $\kappa \le 2^{\aleph_0}$ a parameter for $kap$ is used.   \label{6.2}   \end{theorem}

\noindent {\bf Proof} The case $cf(\kappa) > 2^{\aleph_0}$ follows from Theorems \ref{5.2} and
\ref{5.5}, and the case $cf(\kappa) \le 2^{\aleph_0}$ is covered by the above discussion. As
remarked above, although the parameter $kap$ is needed in the case $\kappa \le 2^{\aleph_0}$,
since `having cardinality $\kappa$' is expressible, we can define when a subset of $IS_n$ is a
possible choice for its interpretation.       $\Box$

\vspace{.1in}    

To complete our analysis of the case $cf(\kappa) \le 2^{\aleph_0}$ we show how ${\cal
M}_{\kappa \lambda \mu}$ is (weakly) interpretable in ${\cal N}_{\kappa \lambda \mu}^2$ in this case
in the sense of Theorem \ref{5.9}. This will suffice to show that the structures ${\cal N}_{\kappa
\lambda \mu}^2$ completely capture the first order theory of the groups
$S_\lambda(\mu)/S_\kappa(\mu)$, which is our goal. Here we use a modification of the definition of a
$k$-representation of a tuple $({\overline t}, {\overline \beta}, {\overline h})$ in ${\cal
M}_{\kappa \lambda \mu}$. Recall that without the assumption $cf(\kappa) > 2^{\aleph_0}$ we
only know that $\overline h$ is a tuple of ${\cal E}_n$-classes of functions, which is one reason
for the altered definition. Another point is that we need to capture the eventual behaviour of two
well-order-types, namely $Card$ above {\em and}  below $\kappa$. If $\kappa = \aleph_\beta$ and
$\gamma$ is least such that $\beta = \gamma + \alpha^*$, we let $Card^* = \{ \nu \in Card: \nu = 0
\vee \aleph_\gamma \le \nu < \lambda \}.$ Then a $k$-representation of $({\overline t}, {\overline
\beta}, {\overline h})$ is defined to be any tuple of the form $(A, <^*, {\overline g}, {\overline t},
{\overline b}, {\overline H})$ such that

$<^*$ well-orders $IS_2$,

$A \subseteq IS_2, b_i \in A, c_0, c_\kappa \in A$,

if $h_i \in F_j$ then $H_i: IS_j \rightarrow A$,

$\overline g$ is a tuple of the form $(g_0, \ldots,g_{k-1},g^0,\ldots, g^{k-1}, g_0^*,
\ldots,g_{k-1}^*,g^{*0},\ldots, g^{* \, k-1})$ where $g_i, g^i: A' \cup \{ \infty \}
\rightarrow IS_2, g_i^*, g^{* \, i}: A'' \rightarrow IS_2$, where $A' = \{ a \in A: a \le^*
c_\kappa\}$ and $A'' = \{ a \in A: c_\kappa < a \}$,

and for some 1--1 order-preserving map $\theta: A \rightarrow Card^*$, 

$\theta$ takes $c_0$ to 0, $c_\kappa$ to $\kappa$, and $b_i$ to $\beta_i$ for each $i$,

if $h_i \in F_j$ then for some $h_i': IS_j \rightarrow Card^*, (h_i')_{{\cal E}_j} = h_i$ and
$(\forall t \in IS_j)\theta(H_i(t))$ $= h_i'(t)$,

the order-types of $\{t \in IS_2: t <^* g_j(a)\}$ and $\{t \in IS_2: t <^* g^j(a)\}$ are equal to
$\gamma (\theta (a), \theta (A'))_{[j]}$ and $\gamma (\theta (a), \theta (A'))^{[j]}$ respectively,
for each $a \in A'$,

and the order-types of $\{t \in IS_2: t <^* g_j^*(a)\}$ and $\{t \in IS_2: t <^* g^{*j}(a)\}$ are
equal to $\gamma (\theta (a), \theta (A''))_{[j]}$ and $\gamma (\theta (a), \theta (A''))^{[j]}$
respectively, for each $a \in A''$. 

\begin{theorem} For every (first order) formula $\varphi(x_0, \ldots , x_{n-1})$ of the
language of ${\cal M}_{\kappa \lambda \mu}$ there is an effectively determined integer $k$ and
(second order) formula $\psi(y_0, \ldots, y_{4k+n+1})$ of the language of ${\cal N}^2_{\kappa \lambda
\mu}$ such that for all $\kappa, \lambda,\mu$ with $cf(\kappa) \le 2^{\aleph_0}$, for every
$a_0,\ldots, a_{n-1}$ in ${\cal M}^*_{\kappa \lambda \mu}$ (having the correct sorts) and every 
$k$-representation $c_0, \ldots , c_{4k+n+1}$ of $\overline a$ in ${\cal N}^2_{\kappa \lambda \mu}$,
$${\cal M}^*_{\kappa \lambda \mu} \models \varphi[{\overline a}] \Leftrightarrow {\cal N}^2_{\kappa
\lambda \mu} \models \psi[{\overline c}].$$            \label{6.3}      \end{theorem}

\noindent{\bf Proof} $\;$ We have to indicate the appropriate modifications in the proof of Theorem
\ref{5.9}. We first remark on the analogues of Lemmas \ref{5.7} and \ref{5.8}, which are required
here too. Finding a formula to express the existence of a $k$-representation is much as before. Some
modification is needed in Lemma \ref{5.8}, since we have to allow for the possibility that the
$H_i$ may be ${\cal E}_j$-equivalent, so that the lack of an order-isomorphism between the
corresponding $A$s need not determine whether or not the $k$-representations are isomorphic. This
is handled using an additional existential quantifier.

Proceeding to the main proof, since we now have to work with ${\cal M}_{\kappa \lambda \mu}$ rather
than ${\cal M}^*_{\kappa \lambda \mu}$, there are some extra atomic cases in the induction to
consider. We concentrate on the formula $Eq(x_0)$, as this serves to illustrate the idea. 

Since the structure ${\cal N}^2_{\kappa \lambda \mu}$, and the notion of `$k$-representation', is
different in the cases $\kappa > 2^{\aleph_0}$ and $\kappa \le 2^{\aleph_0}$, we treat the two
separately, starting with the former, in which $cf(\kappa) \le 2^{\aleph_0} < \kappa$. Let $(A,<^*,H)$
be a 0-representation of $h$. Then
$${\cal M}_{\kappa \lambda \mu} \models Eq[h] \Leftrightarrow \sum \{ h'(t): t \in IS_2 \wedge
(t = ((B,g_1,g_2))_{\cong} \rightarrow g_1 \neq g_2) \} < \kappa,$$
where $(h')_{{\cal E}_2} = h$ is as in the definition of $k$-representation, corresponding to $H$,
and this is equivalent to
$$(\exists \nu < \kappa)(\forall t \in IS_2)(t = ((B,g_1,g_2))_{\cong} \rightarrow h'(t) \le \nu).$$
For if $(\forall t \in IS_2)(t = ((B,g_1,g_2))_{\cong} \rightarrow h'(t) \le \nu)$ then $\sum
\{h'(t): t \in IS_2 \wedge (t = ((B,g_1,g_2))_{\cong} \rightarrow g_1 \neq g_2)\} \le 2^{\aleph_0}.\nu
< \kappa$ as $2^{\aleph_0}, \nu < \kappa$. And if $\{h'(t): t \in IS_2 \wedge (t =
((B,g_1,g_2))_{\cong} \rightarrow g_1 \neq g_2)\}$ is unbounded in $Card_{< \kappa}$, then $\sum
\{h'(t): t \in IS_2 \wedge (t = ((B,g_1,g_2))_{\cong} \rightarrow g_1 \neq g_2)\} \ge \sup \{h'(t): t
\in IS_2 \wedge (t = ((B,g_1,g_2))_{\cong} \rightarrow g_1 \neq g_2)\} \ge \kappa$.

Therefore
\begin{eqnarray*} 
{\cal M}_{\kappa \lambda \mu} \models Eq[h] & \Leftrightarrow & (\exists \nu < \kappa)(\forall t \in
IS_2)(Eq^1(t) \rightarrow h'(t) \le \nu) \\    
& \Leftrightarrow & {\cal N}^2_{\kappa \lambda \mu} \models (\exists y \in IS_2)(\forall z \in
IS_2)(y \in A \wedge y <^* {\underline c}_\kappa \wedge (Eq^1(z)  \\
& & \hspace{2in}  \rightarrow H(t) \le^* y)),
\end{eqnarray*} 
and this provides the desired formula $\psi(y_0, y_1, y_2)$.

Now turning to the case where $\kappa \le 2^{\aleph_0}$ we find that
\begin{eqnarray*} 
{\cal M}_{\kappa \lambda \mu} \models Eq[h] & \Leftrightarrow & (\exists \nu < \kappa)(\forall t \in
IS_2)(t = ((B,g_1,g_2))_{\cong} \rightarrow h'(t) \le \nu) \\    
                &  & \wedge |\{t \in IS_2: h'(t) \neq 0\} < \kappa  \\
                &  & ( \wedge(\forall t \in IS_2 - fin)(h'(t) = 0) \mbox{ when } \kappa = \aleph_0).
\end{eqnarray*} 
The second clause can be expressed by using $kap$; one says that there is a 1--1 function from $\{t:
H(t) \neq c_0\}$ into $kap$, and that no such function is onto. Similarly, when $\kappa = \aleph_0$,
the final clause is expressed by $(\forall t \in IS_2 - fin)(H(t) \neq c_0)$. 

In conclusion we note that although $kap$ is used here, for any two possible choices for it, we can
define when the representation of some object (for instance an ordinal) under the two values really
represents the same object, and so the apparent arbitrariness is inessential.     $\Box$

\vspace{.1in}

\noindent {\bf Further remarks}

We first remark here that if $\kappa$ is a successor cardinal, then the analysis at once becomes much
easier. For if $cf(\kappa) > 2^{\aleph_0}$ then we may apply the results of section 5, and if
$cf(\kappa) \le 2^{\aleph_0}$ then $\alpha^* = 1$ and the extra sorts of the structure ${\cal
N}^2_{\kappa \lambda \mu}$ play no essential part. Note however that although we can distinguish
these two cases ($cf(\kappa) > 2^{\aleph_0}, cf(\kappa) \le 2^{\aleph_0}$), we cannot distinguish
when $\kappa$ is a successor. For as remarked at the beginning of section 5, if $cf(\kappa_1), 
cf(\kappa_2) \ge (2^{\aleph_0})^+$ then ${\cal M}_{\kappa_1 \kappa_1^+ \mu_1} \cong {\cal M}_{\kappa_2
\kappa_2^+ \mu_2}$, but $cf(\kappa) > 2^{\aleph_0}$ is compatible both with $\kappa$ a successor and
$\kappa$ singular.

Arising out of this, we further note that in the general case, (if $\lambda > \kappa^+$), ${\cal
M}_{\kappa \lambda \mu}$ and the disjoint sum of ${\cal M}_{\kappa^+ \lambda \mu}$ and ${\cal
M}_{\kappa \kappa^+ \mu}$ are  bi-interpretable, and so we can separate our problem into two parts,
the first as in the previous paragraph, and the second of which is the true content of section 6.

\vspace{.1in}

\noindent {\bf Conclusions}

In studying the elementary theory of the groups $G = S_\lambda(\mu)/S_\kappa(\mu)$ where $\aleph_0
\le \kappa < \lambda \le \mu^+$ we  distinguish the following eight cases (by first order sentences of
the language of group theory):

First we distinguish the cases $\lambda \le \mu$ and $\lambda = \mu^+$. In each of these, the
elementary theory of $G$ is determined just by the values of $\kappa$ and $\lambda$. Then we
consider the cases
$$cf(\kappa) > 2^{\aleph_0}, \;\; cf(\kappa) \le 2^{\aleph_0} < \kappa, \;\; \aleph_0 < \kappa \le
2^{\aleph_0}, \; \mbox{ and } \; \kappa = \aleph_0.$$
In each case we form a many-sorted second order structure ${\cal N}^2_{\kappa \lambda \mu}$ whose
sorts all have cardinality $\le 2^{\aleph_0}$, which captures the first order theory of $G$, meaning
that 
$$S_{\lambda_1}({\mu_1})/S_{\kappa_1}({\mu_1}) \equiv S_{\lambda_2}({\mu_2})/S_{\kappa_2}({\mu_2})
\Leftrightarrow {\cal N}^2_{{\kappa_1} {\lambda_1} {\mu_1}} \equiv {\cal N}^2_{{\kappa_2} {\lambda_2}
{\mu_2}}.$$ 
For $cf(\kappa) > 2^{\aleph_0}$ we just require information about $Card^-$; in the other cases,
information about the (large enough) cardinals below $\kappa$ is also represented, and when
$\kappa \le 2^{\aleph_0}$ we also require extra unary predicate(s) on $IS_2$. We summarize this by
the general form of Corollary \ref{5.10}:

\begin{corollary} If $\kappa_1 < \lambda_1 \le \mu_1^+$ and $\kappa_2 < \lambda_2 \le \mu_2^+$ then
the following are equivalent:

(i) ${\cal M}_{\kappa_1 \lambda_1 \mu_1} \equiv {\cal M}_{\kappa_2 \lambda_2 \mu_2}$, 

(ii) ${\cal N}^2_{\kappa_1 \lambda_1 \mu_1} \equiv {\cal N}^2_{\kappa_2 \lambda_2 \mu_2}$,

(iii) $S_{\lambda_1}(\mu_1)/S_{\kappa_1}(\mu_1) \equiv S_{\lambda_2}(\mu_2)/S_{\kappa_2}(\mu_2)$. 
                                      \label{6.4}      \end{corollary}

For the case where $\lambda \le \mu$ and $cf(\kappa) > 2^{\aleph_0}$ the following holds: For any
given ordinals $\alpha_l, \alpha^l < \Omega$ there is a first order theory $T$ in the language of
group theory such that
$$  \begin{array}{c}
\mbox{if $2^{\aleph_0} < cf(\kappa) \le \mu, \lambda \le \mu, \kappa = \aleph_\beta, \lambda =
\aleph_\gamma$, $\beta + \alpha = \gamma$, and $\alpha_{[n]} = \alpha_n, \alpha^{[n]} = \alpha^n$
}\\
\mbox{for each $n$, then the first order theory of the group $S_\lambda(\mu)/S_\kappa(\mu)$ is equal
to $T$,}
     \end{array}  $$
with similar statements in the other cases (including reference to the $\alpha^*_{[n]},
\alpha^{*[n]}$ and $kap, fin$ as appropriate).

\vspace{.1in}

Finally we remark on quotients by alternating and trivial groups. The class $\{ S_\lambda(\mu)/A(\mu):
\aleph_0 \le \lambda \le \mu^+ \}$ of quotients by alternating groups is definable in the class of
all quotients of symmetric groups, being precisely those with non-trivial centre. Moreover since the
centre of $S_\lambda(\mu)/A(\mu)$ is just $S_\omega(\mu)/A(\mu)$, which has order $2, 
S_\lambda(\mu)/S_\omega(\mu)$ can be easily interpreted in  $S_\lambda(\mu)/A(\mu)$. It follows that
if  $S_{\lambda_1}(\mu_1)/A(\mu_1) \equiv S_{\lambda_2}(\mu_2)/A(\mu_2)$ then 
$S_{\lambda_1}(\mu_1)/S_\omega(\mu_1) \equiv S_{\lambda_2}(\mu_2)/S_\omega(\mu_2)$, but whether the
converse is true is not at present clear, (though, as we have seen, the class $\{
S_\lambda(\mu)/S_\omega(\mu): \lambda, \mu \}$ {\em is} definable in $\{ S_\lambda(\mu)/S_\kappa(\mu):
\kappa, \lambda, \mu \}$). The quotients by trivial groups are just the normal subgroups of
$S_\lambda(\mu)$, which were studied in \cite{Shelah1} and \cite{Shelah2}. These may be distinguished
from the other `genuine' quotient groups we have studied (as in \cite{Truss}) by means of the sentence
$$\exists x(x \neq 1 \wedge x^2 = 1 \wedge (\forall y)((xx^y)^2 = 1 \vee (xx^y)^3 = 1))$$
(which says that there is a transposition).

\end{document}